\documentclass[9pt]{amsart}
\usepackage{amsfonts}
\usepackage{graphicx}
\usepackage[all]{xy}
\usepackage{amssymb,extarrows}
\usepackage{amscd}
\usepackage{amsmath,amsfonts,amssymb,amsthm,mathrsfs,xypic}
\newtheorem{thm}{Theorem}[section]

\newtheorem{cor}[thm]{Corollary}
\newtheorem{lem}[thm]{Lemma}
\newtheorem{exm}[thm]{Example}

\newtheorem{prop}[thm]{Proposition}
\newtheorem{rem}[thm]{Remark}
\newtheorem{defn}[thm]{Definition}
\numberwithin{equation}{section}

\newcommand{\s}{\hfill\blacksquare}
\newcommand{\End}{\operatorname{End}}
\newcommand{\add}{\operatorname{add}}
\newcommand{\Ima}{\operatorname{Im}}
\newcommand{\Hom}{\operatorname{Hom}}
\newcommand{\Ext}{\operatorname{Ext}}

\begin{document}
\title [Categorical resolution]{Categorical resolutions \\ of a class of derived categories}
\author [Pu Zhang]
{Pu Zhang}
\thanks{Supported by NSF China (Grant No. 11271251 and
11431010).}
\thanks {pzhang$\symbol{64}$sjtu.edu.cn}\maketitle
\begin{center}
Department of Mathematics \\ Shanghai Jiao Tong University, Shanghai
200240, China
\end{center}
\begin{abstract} \ By using the relative derived categories,
we prove that if an Artin algebra $A$ has a module $T$ with ${\rm
inj.dim}T<\infty$ such that $^\perp T$ is finite, then the bounded
derived category $D^b(A\mbox{-}{\rm mod})$ admits a categorical
resolution in the sense of [Kuz], and a categorical
desingularization in the sense of [BO]. For CM-finite Gorenstein
algebra, such a categorical resolution is weakly crepant. The
similar results hold also for $D^b(A\mbox{-}{\rm Mod})$.

 \vskip10pt

\end{abstract}

\section {\bf Introduction}

\subsection{} A categorical resolution of an algebraic variety comes from looking for
a minimal resolution of singularities. The functor
$D^b(\widetilde{X})\rightarrow D^b(X)$ induced by a resolution of
singular variety $X$ enjoys some remarkable properties. This
motivates the study of a categorical resolution of a triangulated
category. A. Bondal and D. Orlov [BO, Section 5] define {\it a
categorical desingularization of triangulated category}  $\mathcal
D$ to be a pair $(D^b(\mathcal A), \mathcal K)$, where $\mathcal A$
is an abelian category of finite homological dimension, and
$\mathcal K$ a thick subcategory, such that $D\cong D^b(\mathcal
A)/\mathcal K$. A. Kuznetsov [Kuz] defines {\it a categorical
resolution} of $\mathcal D$ to be a triple $(\widetilde{\mathcal D},
\pi_*, \pi^*),$  where $\widetilde{\mathcal D}$ is an admissible
subcategory of $D^b(\widetilde{X})$ with $\widetilde{X}$ a smooth
variety, $\pi_*: \widetilde{\mathcal D} \rightarrow \mathcal D$ and
$\pi^*: \mathcal D_{\rm perf} \rightarrow \widetilde{\mathcal D}$
are triangle functors satisfying $(ii)$ and $(iii)$ in Definition
\ref{defcatres} below. If $\pi^*$ is right adjoint to $\pi_*$, then
it is called weakly crepant. M. Van den Bergh [Van] defines a
non-commutative crepant resolution, this induces a categorical
desingularization and a weakly crepant categorical resolution. For
some of the other influential works in this area we refer to [Ab],
[BKR], [BLV], [Kal], [Lun], and [SV].

 \vskip5pt

\subsection{} In this paper we combine Kuznetsov's definition with
Bondal-Orlov's one. See Definition \ref{defcatres}. The reasons are:
this is the case for a proper birational resolution of an algebraic
variety of rational singularity, it contains more information, and
applies to our purpose.

\vskip5pt

Two comments on Definition \ref{defcatres} are in order. First, we
usually need to explicitly determine the perfect subcategory
$\mathcal D_{\rm perf}$ for a work. For an abelian category
$\mathcal A$  with enough projective objects, we have $D^b_{\rm
perf}(\mathcal A): = D^b(\mathcal A)_{\rm perf}\supseteq
K^b(\mathcal P),$  and $D^b_{\rm perf}(\mathcal A) = K^b(\mathcal
P)$ in many important cases, where $\mathcal P$ is the full
subcategory of $\mathcal A$ consisting of projective objects.  For
examples, if $\mathcal A$ is the module category of a ring, or the
finitely generated module category of an Artin algebra. Propositions
\ref{finitely filtrated} and \ref{thick} below say that this is also
the case if $\mathcal A$ is finitely filtrated, or $\mathcal A$ is
cocomplete.

\vskip5pt

Second, there are several ways for defining the smoothness of a
triangulated category ([BO], [KS], [Kuz],  [Lun], and [TV]). In this
paper a triangulated category is smooth, if it is
triangle-equivalent to $D^b(\mathcal A)$ with $\mathcal A$ an
abelian category such that $D^b_{sg}(\mathcal A):=D^b(\mathcal
A)/D^b_{\rm perf}(\mathcal A)=0$ (see Definition \ref{smooth}). This
singularity category $D^b_{sg}(\mathcal A)$ is in the sense of
R.O.Buchweitz [Buch] and D.Orlov [O2]. It is invariant under
triangle-equivalences. An algebraic variety $X$ is smooth if and
only if $D^b_{sg}(X) = 0$ ([O2]). Definition \ref{smooth} is almost
same with the one in [BO, Section 5]. The reason for a replacement
of finite homological dimension by $D^b_{sg}(\mathcal A)= 0$ is that
in general the invariance of finite homological dimension under
derived equivalences seems to be not known, although this is true in
many important cases. Definition \ref{smooth} is also a slight
change of [Kuz, Definition 3.1], since in this paper we do not need
to consider admissible subcategories. We hope that such a small
change will not affect the minimality.
\subsection{} Given a full subcategory $\mathcal C$ of an abelian category $\mathcal A$,
the $\mathcal C$-relative derived category $D_\mathcal C^b(\mathcal
A)$ is the Verdier quotient of the  homotopy category $K^b(\mathcal
A)$ by the thick subcategory $K^b_{{\mathcal C} ac}(\mathcal A)$
consisting of bounded $\mathcal C$-acyclic complexes. It is in fact
a kind of the derived category of an exact category in the sense of
A. Neeman [N], and has been studied in different setups, see for
examples A. B. Buan [Bu], [GZ], X.W.Chen [C2], and J.Asadollahi,
R.Hafezi, and R.Vahed [AHV]. If $\mathcal C$ is contravariantly
finite in $\mathcal A$, then $D_\mathcal C^b(\mathcal A)\cong K^{-,
\mathcal C b}(\mathcal C)$. See Proposition \ref{gpder} for the
details. We observe that the derived category $D^b(\mathcal A)$ can
be described via the relative derived category $D_\mathcal
C^b(\mathcal A)$ (see Theorem \ref{relativedescription} below).

\vskip5pt

\noindent {\bf Theorem 1} \ {\it Let $\mathcal A$ be an abelian
category with enough projective objects, and $\mathcal C$ a
resolving contravariantly finite subcategory of $\mathcal A$. Then
we have a triangle-equivalence $D^b(\mathcal A)\cong D_\mathcal
C^b(\mathcal A)/K^b_{ac}(\mathcal C).$}

\vskip5pt

Thus, we get a functor $\pi_*: D_\mathcal C^b(\mathcal A)\rightarrow
D^b(\mathcal A)$, given by the Verdier functor. As we mentioned
before, one can roughly think $D^b_{\rm perf}(\mathcal A)$ as
$K^b(\mathcal P)$. Since $D_\mathcal C^b(\mathcal A)\cong K^{-,
\mathcal C b}(\mathcal C)$ and $K^b(\mathcal P) \subseteq K^{-,
\mathcal C b}(\mathcal C)$, we have another functor $\pi^*:
K^b(\mathcal P) \rightarrow D_\mathcal C^b(\mathcal A)$,  given by
the embedding. So, if $D_\mathcal C^b(\mathcal A)$ is smooth, plus
some other required properties of $\pi_*$ and $\pi_*$, then the
triple $(D_\mathcal C^b(\mathcal A), \pi_*, \pi^*)$ can be served as
a categorical resolution of $D^b(\mathcal A)$ in our consideration.
We will see below that for some kinds of $\mathcal A$ and $\mathcal
C$, this machinery works.

\vskip5pt

\subsection{} Let $A$ be an Artin algebra, $A\mbox{-}{\rm
mod}$ (resp. $A\mbox{-}{\rm Mod}$) the category of left $A$-modules
(resp. finitely generated left $A$-modules). We say that $A$ is {\it
representation-finite}, if $A\mbox{-}{\rm mod}$ has only finite many
pairwise non-isomorphic indecomposable objects. By a theorem of M.
Auslander, in this case the indecomposable objects in $A\mbox{-}{\rm
Mod}$ coincide with the ones indecomposable objects in
$A\mbox{-}{\rm mod}$ ([A2]). For  $T\in A\mbox{-}{\rm mod}$, let
$^\perp T$ denote the full subcategory of $A\mbox{-}{\rm mod}$
consisting of $A$-modules $X$ such that $\Ext^i_A(X, T) = 0, \
\forall \ i\ge 1.$ By ${\rm add} T$ we denote the full subcategory
of $A$-{\rm mod} consisting of the direct summands of finite direct
sums of copies of $T$, and by ${\rm Add} T$ the full subcategory of
$A$-{\rm Mod} consisting of the direct summands of arbitrary direct
sums of copies of $T$. We say that $^\perp T$ is finite, if there
are only finitely many pairwise non-isomorphic indecomposable
$A$-modules in $^\perp T$, or equivalently, there is a module $M$ in
$A$-mod such that $^\perp T = {\rm add} M$. The following result on
the endomorphism algebras of finite global dimension is another key
step for the categorical resolution in this paper.

\vskip5pt

\noindent {\bf Theorem 2} \ Let $A$ be an Artin algebra,  $T$ and
$M$ modules in $A\mbox{-}{\rm mod}$ such that  $^\perp T = {\rm add}
M.$ Put $B: = (\End_A M)^{op}$. Then for each positive integer $r\ge
2$, ${\rm gl.dim} B\le r$ if and only if  ${\rm inj.dim} T\le r$.

\vskip5pt

Two special cases of Theorem 1 are well-known. If $A$ is
representation-finite and $T$ is an injective $A$-module, then the
corresponding $B$ is the Auslander algebra {\rm([ARS])}. If $A$ is
{\rm CM}-finite Gorenstein algebra and  $T =\ _AA$, then the
corresponding $B$ is the relative Auslander algebra of $A$ (see
[LZ], [Bel2], and [Leu]).

\vskip10pt

\subsection{} Assume that $A$ is of infinite global
dimension. If there exist modules $T$ and $M$ in $A$-mod  with ${\rm
inj.dim} T< \infty$,  such that $^\perp T = {\rm add} M$, then
$D^b(B\mbox{-}{\rm mod})$ is smooth by Theorem 2, where $B =
(\End_AM)^{op}$. While $D^b(B\mbox{-}{\rm mod})$ is
triangle-equivalent to the relative derived category $D_{{\rm
add}M}^b(A\mbox{-}{\rm mod})$, so by the comment after Theorem 1 we
have a triple $(D^b(B\mbox{-}{\rm mod}), \pi_*, \pi^*)$. The main
results (Theorems \ref{mainresult}, \ref{finrem} and
\ref{mainresult3}) of this paper say that this triple gives a
categorical resolution of $D^b(A\mbox{-}{\rm mod})$.

\vskip5pt

\noindent{\bf Theorem 3.} \ {\it Let $A$ be an Artin algebra with
${\rm gl.dim}A = \infty$.

\vskip5pt

$(i)$ \ Assume that there are modules $T$ and $M$ in $A\mbox{-}{\rm
mod}$ with ${\rm inj.dim} T< \infty$, such that $\ ^\perp T = {\rm
add}M$. Then $D^b(B\mbox{-}{\rm mod})$ is a categorical resolution
$D^b(A\mbox{-}{\rm mod})$, where $B = (\End_AM)^{op}$.

\vskip5pt

$(ii)$ \  Assume that there are modules $T$ and $M$ in
$A\mbox{-}{\rm mod}$ with ${\rm inj.dim} T< \infty$, such that
$^{\perp_{\rm big}} ({\rm Add}T) = {\rm Add} M$. Then
$D^b(B\mbox{-}{\rm Mod})$ is a categorical resolution
$D^b(A\mbox{-}{\rm Mod})$, where $B = (\End_AM)^{op}$.

\vskip10pt

$(iii)$ If $A$ is a {\rm CM}-finite Gorenstein algebra, and $B$ the
relative Auslander algebra of $A$, then $D^b(B\mbox{-}{\rm mod})$ is
a weakly crepant categorical resolution of $D^b(A\mbox{-}{\rm
mod})$, and $D^b(B\mbox{-}{\rm Mod})$ is a weakly crepant
categorical resolution of $D^b(A\mbox{-}{\rm Mod})$.}

\vskip5pt

We remark that in Theorem 3$(iii)$ if $A$ is in addition a
commutative local ring, then G. J. Leuschke [Leu, Section 3] has
observed a connection with non-commutative crepant resolution in the
sense of M. Van den Bergh [Van].

\subsection {} The paper is organized as
follows. In $\S2$ we recall the main definitions and facts used. For
an abelian category $\mathcal A$ with enough projective objects, in
$\S 3$ we give two frequently used cases of $\mathcal A$, such that
$D^b_{\rm perf}(\mathcal A) = K^b(\mathcal P)$. In $\S 4$ we prove
Theorem 2. In $\S 5$ we recall some points of the relative derived
categories; and in $\S 6$ we prove Theorem 1 and other facts, which
provide the adjointness of the functors appeared in the categorical
resolution. In $\S 7$ we prove the main results with some
consequences.

\vskip5pt

\section{\bf Preliminaries}
\subsection{The key formula} Let $\mathcal A$ be an abelian category. For $*\in\{{\rm blank}, -, b\}$, let $K^*(\mathcal{A})$ and
$D^*(\mathcal{A})$ be the corresponding homotopy category and the
derived category of $\mathcal A$, respectively. For complexes $X$
and $Y$, let $\operatorname{Hom}_\mathcal A(X, Y)$ be the Hom
complex. Then we have the key formula
$\operatorname{Hom}_{K(\mathcal A)}(X, Y[n])$ $ = {\rm
H}^n\operatorname{Hom}_{\mathcal A} (X, Y), \ \forall \ n\in\Bbb Z.$

\vskip5pt

\subsection{Verdier quotients}\ Let $\mathcal B$ be a
triangulated subcategory of triangulated category $\mathcal K$.
Thus, in particular, $\mathcal B$ is a full subcategory of $\mathcal
K$, and closed under isomorphisms ([N]). By definition a morphism
$f: X\longrightarrow Y$ of the Verdier quotient $\mathcal K/\mathcal
B$ is an equivalence class of right fractions $a/s$, where $s: Z
\Longrightarrow X$ and $a: Z\longrightarrow Y$ are morphisms of
$\mathcal K$, such that `` the mapping cone"  of $s$ belongs to
$\mathcal B$. Let $Q: \mathcal K \longrightarrow \mathcal K/\mathcal
B$ be the localization functor sending an object $X$ to $X$ itself,
and sending a morphism $a: X\longrightarrow Y$ to $a/{\rm Id}_X$.
Then $Q$ is a triangle functor with $Q(\mathcal B) = 0$; and if $F:
\mathcal K\longrightarrow \mathcal T$ is a triangle functor with
$G(\mathcal B) = 0$, then there is a unique triangle functor $G:
\mathcal K/\mathcal B\longrightarrow \mathcal T$ such that $F = GQ$.
Thus $Q(X) \cong 0$ if and only if $X$ is a direct summand of an
object in $\mathcal B$. If $\mathcal B$ is {\it thick in $\mathcal
K$} (i.e., $\mathcal B$ is a triangulated subcategory of $\mathcal
K$ which is closed under direct summands), then $Q(X) \cong 0$ if
and only if $X\in \mathcal B$.

\vskip5pt

We also need the left fraction construction of $\mathcal K/\mathcal
B$: a morphism $f: X\longrightarrow Y$ in $\mathcal K/\mathcal B$ is
an equivalence class of left fractions $s\backslash a$, where $a:
X\longrightarrow Z$ and $s: Y \Longrightarrow Z$ are morphisms of
$\mathcal K$, such that `` the mapping cone" of $s$ belongs to
$\mathcal B$. The localization functor sends a morphism $a:
X\longrightarrow Y$ to ${\rm Id}_Y\backslash a$. Then the Verdier
quotient $\mathcal K/\mathcal B$ constructed via right fractions is
isomorphic to the one constructed via left fractions.

\begin{lem} \label{Verdier} {\rm ([Ver, \mbox{Corollary}  4\mbox{-}3])} \ Let $\mathcal{D}_1$ and $\mathcal{D}_2$
be triangulated subcategories of triangulated category
$\mathcal{C}$, and $\mathcal{D}_1$ a subcategory of $\mathcal{D}_2$.
Then $\mathcal{D}_2/\mathcal{D}_1$ is a triangulated subcategory of
$\mathcal{C}/\mathcal{D}_1$, and
$(\mathcal{C}/\mathcal{D}_1)/(\mathcal{D}_2/\mathcal{D}_1)\cong
\mathcal{C}/\mathcal{D}_2$ as triangulated categories.
\end{lem}

\subsection{Perfect objects of a triangulated category}
Let $\mathcal D$ be a triangulated category, with shift functor
denoted by $[1]$. Following [Kuz], an object $P\in \mathcal D$ is
{\it perfect}, provided that for each object $Y$ of $\mathcal D$,
there are only finitely many $i\in \Bbb Z$ such that $\Hom_\mathcal
D(P, Y[i]) \ne 0$. Denote by $\mathcal D_{\rm perf}$ the full
subcategory of $\mathcal D$ consisting of perfect objects, which is
called {\it the perfect subcategory} of $\mathcal D$. Then $\mathcal
D_{\rm perf}$ is a thick subcategory of $\mathcal D$.

\vskip5pt

This definition comes from the intrinsic characterization of a
perfect complex of $D^b(X)$ of coherent sheaves on algebraic variety
$X$ ([O1, Proposition 1.11]). Its advantage is that a
triangle-equivalence $\mathcal D\longrightarrow \mathcal D'$
restricts to a triangle-equivalence $\mathcal D_{\rm
perf}\longrightarrow \mathcal D'_{\rm perf}$, and hence induces a
triangle-equivalence $\mathcal D/\mathcal D_{\rm
perf}\longrightarrow \mathcal D'/\mathcal D'_{\rm perf}$. For an
abelian category $\mathcal A$ with enough projective objects, one
has $D^b_{\rm perf}(\mathcal A): =  D^b(\mathcal A)_{\rm
perf}\supseteq K^b(\mathcal P),$  and $D^b_{\rm perf}(\mathcal A) =
K^b(\mathcal P)$ in many important cases, where $\mathcal P$ is the
full subcategory of $\mathcal A$ consisting of projective objects.
For examples, this is the case when $\mathcal A = R$-Mod for ring
$R$, or $\mathcal A = A$-mod for an Artin algebra $A$. For the
details see $\S 3$. However, in general we do not know whether
$D^b_{\rm perf}(\mathcal A) = K^b(\mathcal P)$.

\vskip5pt

\subsection{Smooth triangulated categories}
Let $X$ be an algebraic variety. D. Orlov called the Verdier
quotient $D^b_{sg}(X): = D^b(X)/D^b_{\rm perf}(X)$ the singularity
category of $X$, where $D^b_{\rm perf}(X) : = (D^b(X))_{\rm perf}$.
Then $X$ is smooth if and only if $D^b_{sg}(X) = 0$.

\vskip5pt

For any abelian category $\mathcal A$, this kind of Verdier quotient
was introduced by R. O. Buchweitz [Buch, 1.2.2], under the name {\it
the stabilized derived category} of $\mathcal A$. For short,
following [O2], in this paper we call $D^b_{sg}(\mathcal A): =
D^b(\mathcal A)/D^b_{\rm perf}(\mathcal A)$ {\it the singularity
category} of $\mathcal A$. If $\mathcal A$ has enough projective
objects and $D^b_{\rm perf}(\mathcal A) = K^b(\mathcal P)$, then
$D^b_{sg}(\mathcal A) = 0$ if and only if each object of $\mathcal
A$ has a finite projective dimension. Thus, if $\mathcal A = A$-mod
for an Artin algebra $A$, then $D^b_{sg}(\mathcal A) = 0$ if and
only if the global dimension of $A$ is finite.

\vskip5pt

\begin{defn} \label{smooth}  \ A triangulated category
$\mathcal D$ is smooth, if it is triangle-equivalent to
$D^b(\mathcal A)$ with $D^b_{\rm sg}(\mathcal A) = 0$, where
$\mathcal A$ is an abelian category.
\end{defn}
Thus, by the remark at the previous subsection, if a triangulated
category $\mathcal D$ is smooth, then $\mathcal D = \mathcal D_{\rm
perf}$. In particular, $D^b(R\mbox{-}{\rm Mod})$ is smooth if and
only if each $R$-module has a finite projective dimension.

\vskip5pt

Definition \ref{smooth} is almost same with a smooth triangulated
category in the sense of A. Bondal and D. Orlov [BO, Section 5],
where it is triangle-equivalent to $D^b(\mathcal A)$ with $\mathcal
A$ an abelian category of finite homological dimension (i.e., for
each object $X\in\mathcal A$, there are only finite many integers
$i$ such that $\Ext^i_\mathcal A(X, -)\ne 0$). In fact, if $\mathcal
A$ has enough projective objects, then a smooth triangulated
category in the sense of [BO] is smooth in the sense of Definition
\ref{smooth}; and the converse is also true  if $D^b_{\rm
perf}(\mathcal A) = K^b(\mathcal P)$. The reason to make such a
minor change is that  we do not know in general the invariance of
finite homological dimension under derived equivalences, although
this is true in many important cases.

\vskip5pt

Definition \ref{smooth} is also a slight modification of [Kuz,
Definition 3.1], where a smooth triangulated category is
triangle-equivalent to an admissible subcategory of $D^b(X)$ with
$D^b_{sg}(X) = 0$, where $X$ is an algebraic variety. For our
purpose we do not need to consider admissible subcategories.

\subsection{Categorical resolution of a triangulated category}  For a triangle functor $F: \mathcal T \longrightarrow \mathcal T'$,
let ${\rm Ker} F$ denote the full subcategory of $\mathcal T$
consisting of objects $K$ with $F(K)\cong 0$. The following
definition is due to A. Bondal and D. Orlov [BO, Section 5] and A.
Kuznetsov [Kuz, Definition 3.2].

\begin{defn} \label{defcatres} \ A categorical resolution of a non-smooth triangulated
category $\mathcal D$ is a smooth triangulated category
$\widetilde{\mathcal D}$, or more precisely, a triple
$(\widetilde{\mathcal D}, \pi_*, \pi^*)$, where $\pi_*:
\widetilde{\mathcal D} \longrightarrow \mathcal D$ and $\pi^*:
\mathcal D_{\rm perf}  \longrightarrow \widetilde{\mathcal D}$ are
triangle functors, such that

\vskip5pt

$(i)$ \ $\pi_*$ induces a triangle-equivalence $\widetilde{\mathcal
D}/{\rm Ker}\pi_*\cong \mathcal D;$

\vskip5pt

$(ii)$ \ $\pi^*$ is left adjoint to $\pi_*$ on $\mathcal D_{\rm
perf}$, that is, there is a functorial isomorphism $\eta_{P, X}:
\Hom_{\widetilde{\mathcal D}}(\pi^*P, X)\cong \Hom_\mathcal D(P,
\pi_*X), \ \ \forall \ P\in \mathcal D_{\rm perf}, \ \forall \ X\in
\widetilde{\mathcal D};$

\vskip5pt

$(iii)$ \ The unit $\eta = (\eta_P)_{P\in \mathcal D_{\rm perf}}:
{\rm Id}_{\mathcal D_{\rm perf}} \longrightarrow \pi_*\pi^*$ is a
natural isomorphism of functors, where $\eta_P$ is the morphism
$\eta_{P, \pi^*P}({\rm Id}_{\pi^*P}): P\longrightarrow \pi_*\pi^*P$
in $\mathcal D$.
\end{defn}

Note that $(ii)$ implies that $\pi^*: \mathcal D_{\rm perf}
\longrightarrow \widetilde{\mathcal D}$ is fully faithful.

\vskip5pt

If $\pi_*: \widetilde{\mathcal D} \longrightarrow \mathcal D$ is
full and dense, then $(i)$ in Definition \ref{defcatres} holds
automatically. However $\pi_*$ usually can not be full.

\vskip5pt

It is well-known that for a complex singular variety $X$ there is a
proper birational resolution of singularities
$\widetilde{X}\longrightarrow X$; and that if
$\widetilde{X}\longrightarrow X$ is a proper birational resolution
of algebraic variety $X$ of rational singularity, then
$D^b(\widetilde{X})$ is a categorical resolution of $D^b(X)$ in the
sense of Definition \ref{defcatres}.

\vskip5pt

\begin{defn} \label{wc} {\rm ([Kuz], Definition
3.4)} \ A categorical resolution $(\widetilde{\mathcal D}, \pi_*,
\pi^*)$ of a triangulated category $\mathcal D$ is weakly crepant if
$\pi^*$ is right adjoint to $\pi_*$ on $\mathcal D_{\rm perf}$, that
is, there is a functorial isomorphism $\Hom_{\widetilde{\mathcal
D}}(X, \pi^*P)\cong \Hom_\mathcal D(\pi_*X, P), \ \ \forall \ P\in
\mathcal D_{\rm perf}, \ \forall \ X\in \widetilde{\mathcal D}.$
\end{defn}

A non-commutative crepant resolution ([Van]) induces a weakly
crepant categorical resolution of a triangulated category.

\subsection{Gorenstein-projective objects} Let $\mathcal A$ be an abelian category
with enough projective objects, and $\mathcal P = \mathcal
P(\mathcal A)$ the full subcategory of $\mathcal A$ consisting of
projective objects. {\it A complete $\mathcal A$-projective
resolution} is an exact sequence ${P}^{\bullet}=\cdots
\longrightarrow P^{-1}\stackrel{d^{-1}}\longrightarrow P^0
\stackrel{d^0}\longrightarrow P^{1}\longrightarrow \cdots$ with each
$P^i\in\mathcal P$, such that ${\rm Hom}_\mathcal A({P}^{\bullet},
P)$ stays exact for each $P\in \mathcal P$. An object $G$ of
$\mathcal A$ is {\it Gorenstein-projective} if there is a complete
$\mathcal A$-projective resolution ${P}^{\bullet}$ such that $G\cong
{\rm Im} d^0$ (E. E. Enochs and O. M. G. Jenda [EJ]). Denote by
$\mathcal{GP}(\mathcal A)$ the full subcategory of $\mathcal A$
consisting of Gorenstein-projective objects.

\vskip5pt

A full subcategory $\mathcal X$ of $\mathcal A$ is {\it resolving}
([AB]), provided that $\mathcal X\supseteq \mathcal P$, $\mathcal X$
is closed under extensions and direct summands, and that $\mathcal
X$ is closed under the kernels of epimorphisms. A resolving
subcategory is of course additive. Then $\mathcal{GP}(\mathcal A)$
is resolving; and $\mathcal{GP}(\mathcal A)$ is closed under
arbitrary direct sums if $\mathcal A$ is {\it cocomplete}, i.e.,
$\mathcal A$ has arbitrary direct sums (see [AR] and [Hol]).

\vskip5pt

{\it A Frobenius category} $\mathcal B$ is an exact category ([Q,
$\S$2]) with enough projective objects and enough injective objects,
such that  an object is projective if and only if it is injective
(see [K1]). An important feature is that $\mathcal{GP}(\mathcal A)$
is a Frobenius category, where the projective-injective objects of
$\mathcal{GP}(\mathcal A)$ are exactly the projective objects of
$\mathcal A$ (see [Bel1]). Thus the stable category
$\underline{\mathcal {GP}(\mathcal A)}$ of $\mathcal {GP}(\mathcal
A)$ modulo $\mathcal {P}$ is triangulated ([Hap, p.16]).

\vskip5pt

Recall that an Artin algebra $A$ is {\it CM-finite}, if
$\mathcal{GP}(A\mbox{-}{\rm mod})$ has only finitely many pairwise
non-isomorphic indecomposable objects; and  that $A$ is {\it
Gorenstein}, if ${\rm inj.dim}\ _AA < \infty$ and ${\rm inj.dim} A_A
< \infty$. If this is the case, then ${\rm inj.dim}\ _AA = {\rm
inj.dim} A_A$ ([I]). For a Gorenstein algebra $A$, we have $\mathcal
{GP}(A\mbox{-}{\rm mod}) = \ ^\perp (_AA)$ (see [EJ, Corollary
11.5.3]; or [Z, Lemma 2.4$(iii)$] for a short argument). If $A$ is a
{\rm CM}-finite Gorenstein algebra, then the indecomposable objects
of $\mathcal{GP}(A\mbox{-}{\rm Mod})$ coincide with the
indecomposable objects of $\mathcal{GP}(A\mbox{-}{\rm mod})$, by X.
W. Chen ([C1]).
\subsection{Contravariantly finite subcategories} Let $\mathcal{B}$ be an
additive category, $\mathcal{C}$ a full additive subcategory of
$\mathcal B$, and $X\in\mathcal B$. A morphism $f: C \longrightarrow
X$ with $C\in \mathcal{C}$ is {\it a right
$\mathcal{C}$-approximation} of $X$, if $\Hom_{\mathcal B}(C', f):
\Hom_{\mathcal B}(C', C) \longrightarrow \Hom_{\mathcal B}(C', X)$
is surjective for each $C' \in \mathcal{C}$.  If each object $X\in
\mathcal B$ admits a right $\mathcal{C}$-approximation, then
$\mathcal{C}$ is said to be {\it contravariantly finite in $\mathcal
B$} ([AR]).
\begin{exm}\label{exm} Recall the contravariantly finite
subcategories used in this paper.

\vskip5pt

$(i)$ Let $A$ be an Artin algebra and $M\in A$-mod. Then ${\rm add}
M$ is contravariantly finite in $A$-{\rm mod}; and ${\rm Add} M$ is
contravariantly finite in $A$-{\rm Mod} {\rm(}here $M$ is also
assumed to be finitely generated{\rm)}.

\vskip5pt

$(ii)$ \  If each object of $\mathcal{A}$ has a finite
Gorenstein-projective dimension, then $\mathcal {GP}(\mathcal A)$ is
contravariantly finite in $\mathcal{A}$ {\rm(}{\rm [EJ, Theorem
11.5.1]}, or {\rm [Hol, Theorem 2.10]}{\rm)}.

\vskip5pt

$(iii)$ For an Artin algebra $A$, $\mathcal {GP}(A\mbox{-}{\rm
Mod})$ is contravariantly finite in $A$-{\rm Mod} {\rm(}{\rm [Bel1,
Theorem 3.5]}{\rm)}.

\vskip5pt

$(iv)$ An Artin algebra $A$ is Gorenstein if and only if each
$A$-module has a finite Gorenstein-projective dimension in $A$-{\rm
Mod}, also if and only if each finitely generated $A$-module has a
finite Gorenstein-projective dimension in $A$-{\rm mod} {\rm
([Hos])}. Thus, if $A$ is Gorenstein, then $\mathcal
{GP}(A\mbox{-}{\rm mod})$ is contravariantly finite in $A$-{\rm
mod}.

\vskip5pt

$(v)$ {\rm A. Beligiannis} {\rm [Bel1]} introduced virtually
Gorenstein algebras. A Gorenstein algebra is virtually Gorenstein,
but the converse is not true. However, for a virtually Gorenstein
algebra $A$, $\mathcal {GP}(A\mbox{-}{\rm mod})$ is contravariantly
finite in $A$-{\rm mod} {\rm ([Bel1, Theorem 8.2${\rm (ix)}$])}.

\vskip5pt

$(vi)$  For examples of \ {\rm CM}-finite non-Gorenstein algebras we
refer to {\rm [Rin]}. For a {\rm CM}-finite algebra $A$, $\mathcal
{GP}(A\mbox{-}{\rm mod})$ is contravariantly finite in $A$-{\rm
mod}.
\end{exm}

\section{\bf Perfect subcategory of a triangulated
category}

Throughout this section, $\mathcal A$ is an abelian category with
enough projective objects, $\mathcal P = \mathcal P(\mathcal A)$ the
full subcategory of $\mathcal A$ consisting of projective objects.
We give two classes of $\mathcal A$, such that $D^b_{\rm
perf}(\mathcal A) = K^b(\mathcal P).$
\subsection{} The following characterization of objects in $K^b(\mathcal
P)$ is due to Buchweitz. It also implies that $K^b(\mathcal P)$ is
thick in $D^b(\mathcal A)$.

\begin{lem}\label{buch} {\rm ([Buch, Lemma 1.2.1])} \ Let $\mathcal A$ be an abelian category with enough projective
objects, and $P\in D^b(\mathcal A)$. Then the following are
equivalent

\vskip5pt

$(i)$ \ $P\in K^b(\mathcal P);$

\vskip5pt

$(ii)$ \ there is an integer $i(P)$ such that $\Hom_{D^b(\mathcal
A)}(P, M[i]) = 0$ for each $i\ge i(P)$ and for each object $M$ of
$\mathcal A;$

\vskip5pt

$(iii)$ \ there is a finite subset $I(P)\subseteq \Bbb Z$, such that
$\Hom_{D^b(\mathcal A)}(P, M[j]) = 0$ for each $j\notin I(P)$ and
for each object $M$ of $\mathcal A$.
\end{lem}
\noindent{\bf Proof.} For convenience of the reader, we include an
argument for $(iii) \Longrightarrow (i)$. Let $Q \longrightarrow P$
be a quasi-isomorphism  with $Q\in K^{-, b}(\mathcal P)$. Then there
is an integer $N$ such that ${\rm H}^n Q = 0$  for all $n \le N$. We
claim that there exists an integer $n$ with $n\le N$ such that ${\rm
Im}d^{n}_Q\in \mathcal P$. If this claim is true, then there is a
quasi-isomorphism
\[\xymatrix{Q\ar[d]_{f}:
& \cdots \ar [r] & Q^{n-1} \ar [r] \ar[d]& Q^{n} \ar[d] \ar[r]
 & Q^{n+1}\ar[r]\ar@{=}[d]& \cdots
\\
\tau_{\ge n+1}Q: & \cdots  \ar [r] & 0 \ar [r] &
\operatorname{Im}d^{n}_Q \ar[r]  & Q^{n+1}\ar[r]&
 \cdots.}\] Therefore we have isomorphisms
$P\cong Q \cong \tau_{\ge n+1}Q\in K^b(\mathcal P)$ in $D^-(\mathcal
A)$, and hence we have an isomorphism $P\cong \tau_{\ge n+1} Q \in
K^b(\mathcal P)$ in $D^b(\mathcal A)$.

\vskip5pt

Assume that the claim is not true. Then there exists $-n\notin I(P)$
such that $M:= {\rm
 Im}d^{n}_Q\notin \mathcal P$. Denote by
$\widetilde{d^{n}}: \ Q^{n}\longrightarrow M={\rm
 Im}d^{n}_Q$ the epimorphism induced by $d^n_Q$. Then
$$\widetilde{d^{n}}\in {\rm Ker} ({\rm Hom}_{\mathcal A}(Q^{n}, M)
\stackrel {{\rm Hom}_{\mathcal A}(d^{n-1}, M)}\longrightarrow {\rm
Hom}_{\mathcal A}(Q^{n-1}, M));$$ but
$$\widetilde{d^{n}}\notin {\rm Im} ({\rm Hom}_{\mathcal A}(Q^{n+1}, M) \stackrel {{\rm
Hom}_{\mathcal A}(d^{n}, M)}\longrightarrow {\rm Hom}_{\mathcal
A}(Q^{n}, M))$$ (otherwise ${\rm Im}d^{n}\hookrightarrow Q^{n+1}$
splits, which contradicts  the assumption $M\notin \mathcal P$).
This implies ${\rm H}^{-n}{\rm Hom}_{\mathcal A}(Q, M)\ne 0$.
Therefore
\begin{align*}{\rm Hom}_{D^b(\mathcal A)}(P, M[-n])&\cong
{\rm Hom}_{D^-(\mathcal A)}(Q, M[-n])\\
& \cong {\rm Hom}_{K^-(\mathcal A)}(Q, M[-n])\\& = {\rm H}^{-n}{\rm
Hom}_{\mathcal A}(Q, M) \ne 0.\end{align*} This contradicts the
assumption $(iii)$. \hfill $\s$

\subsection{} We
say that $\mathcal A$ is {\it a finitely filtrated category}, if
there exists finitely many objects $S_1, \cdots, S_m$, such that for
any non-zero object $X$ of $\mathcal A$, there exists a sequence of
monomorphisms
$$0= X_0 \stackrel {f_0} \longrightarrow \cdots
\longrightarrow X_{n-1}\stackrel {f_{n-1}} \longrightarrow X_n = X$$
such that ${\rm Coker} f_i \in\{S_1, \cdots, S_m\}, \ 0\le i\le
n-1$.

\vskip5pt

For example, for an Artin algebra $A$, $A$-mod is finitely
filtrated.

\vskip5pt

\begin{prop} \label{finitely filtrated}  \ Let $\mathcal A$ be a finitely
filtrated category. Then  $D^b_{\rm perf}(\mathcal A) = K^b(\mathcal
P)$.
\end{prop}
\noindent{\bf Proof.} \ We only justify $D^b_{\rm perf}(\mathcal
A)\subseteq K^b(\mathcal P)$. Let $P\in D^b_{\rm perf}(\mathcal A)$.
By assumption $\mathcal A$ is finitely filtrated by some objects
$S_1, \cdots, S_m$. Put $S: = S_1\oplus \cdots \oplus S_m$. Since
$P\in D^b_{\rm perf}(\mathcal A)$, there are only finitely many
$i\in \Bbb Z$ such that $\Hom_{D^b(\mathcal A)}(P, S[i]) \ne 0$.
Denote by $I(P)$ the finite set of such integers $i$'s. Then
$\Hom_{D^b(\mathcal A)}(P, S[j]) = 0$ for $j\notin I(P)$. Since each
object $M\in\mathcal A$ has a filtration with factors belonging to
$\{S_1, \cdots, S_m\}$, and since each short exact sequence in
$\mathcal A$ gives rise a distinguished triangle in $D^b(\mathcal
A)$, it follows that $\Hom_{D^b(\mathcal A)}(P, M[j]) = 0$ for each
$j\notin I(P)$ and for each object $M$ of $\mathcal A$. Thus  $P\in
K^b(\mathcal P)$ by Lemma \ref{buch}.  \hfill $\s$

\vskip10pt

\subsection{} If $\mathcal A$
 has arbitrary direct sums, then we have the same conclusion as in Proposition \ref{finitely filtrated}.
 It in particular say $D^b_{\rm perf}(R\mbox{-}{\rm Mod}) = K^b(\mathcal
P(R\mbox{-}{\rm Mod}))$, where $R$ is a ring.

\vskip10pt

\begin{prop} \label{thick} \ Let $\mathcal A$ be a cocomplete abelian
category with enough projective objects. Then $D^b_{\rm
perf}(\mathcal A) = K^b(\mathcal P)$.
\end{prop}

\noindent{\bf Proof.} \  We only prove  $D^b_{\rm perf}(\mathcal
A)\subseteq K^b(\mathcal P)$. The idea of the proof could be found
from J. Rickard [Ric, Proposition 6.2]. Let $P\in D^b_{\rm
perf}(\mathcal A)$. Take a quasi-isomorphism $Q \longrightarrow P$
with $Q\in K^{-, b}(\mathcal P)$. Thus, there is an $N\in\Bbb Z$
such that ${\rm H}^nP = 0, \ \forall \ n\le N$. As in the proof of
Proposition \ref{buch} it suffices to prove that there exists an
integer $n$ with $n\le N$ such that ${\rm Im}d^{n}_Q\in \mathcal P$.

\vskip5pt

Otherwise, ${\rm
 Im}d^n_Q\notin \mathcal P$ for each $n\le N$.  Since $\mathcal A$ has infinite direct
sums, we could put $M:= \bigoplus\limits_{n\le N}{\rm
 Im}d_Q^{n}\in \mathcal A$. Since ${\rm Im}d_Q^n\ne 0$, we have a
non-zero epimorphism $\widetilde{d^n}: Q^{n}\longrightarrow {\rm
Im}d_Q^{n}$, which induces a non-zero morphism
$$f: Q^{n}\longrightarrow M = \bigoplus\limits_{j\le N}{\rm
 Im}d_Q^{j} = {\rm Im}d_Q^{n}\oplus (\bigoplus\limits_{j\le N, j\ne n}{\rm
 Im}d_Q^{j}).$$ Clearly $f$ induces a chain map $Q
\longrightarrow M[-n]$.  Since ${\rm Im}d_Q^n\notin \mathcal P$, it
follows that this chain map is not null homotopic. This shows ${\rm
Hom}_{K^-(\mathcal A)}(Q, M[-n])\ne 0$ for each integer $n$ with
$n\le N$, and hence
\begin{align*}{\rm Hom}_{D^b(\mathcal A)}(P, M[-n])&\cong
{\rm Hom}_{D^-(\mathcal A)}(Q, M[-n])\\
& \cong {\rm Hom}_{K^-(\mathcal A)}(Q, M[-n])\ne 0.\end{align*} In
other words, we get infinitely many integers $i$ such that ${\rm
Hom}_{D^b(\mathcal A)}(P, M[i])\ne 0$. This contradicts the
assumption $P\in D^b_{\rm perf}(\mathcal A)$. \hfill $\s$

\vskip 10pt

\section{\bf Global dimension of a class of endomorphism algebras}

\subsection {} Let $\mathcal A$ be an abelian category with enough
projective objects, and $X$ an object of $\mathcal A$. The global
dimension ${\rm gl.dim}\mathcal A$ is the supreme of the projective
dimension ${\rm proj.dim} X$, where $X$ runs over all the objects of
$\mathcal A$. For a ring $R$, ${\rm gl.dim} (R\mbox{-}{\rm Mod})$ is
exactly the supreme of ${\rm proj.dim} M$, where $M$ runs over all
the cyclic left $R$-modules (see [A, Theorem 1]). Thus, if $R$ is
left noetherian, then ${\rm gl.dim}(R\mbox{-}{\rm Mod}) = {\rm
gl.dim} (R\mbox{-}{\rm mod}),$ which will be denoted by ${\rm
gl.dim}R$. Thus, for Artin algebra $A$,  ${\rm gl.dim}A$ is just the
maximum of ${\rm proj.dim} S(i), \ 1\le i\le n$, where $\{S(1),
\cdots, S(n)\}$ is a complete set of pairwise non-isomorphic simple
$A$-modules.

\vskip10pt

\subsection{} \ Let $A$ be an Artin algebra, and $M\in
A\mbox{-}{\rm mod}$. The  functor $\Hom_{A}(M, -): A\mbox{-}{\rm
mod} \longrightarrow B\mbox{-}{\rm mod}$ induces an equivalence
between $\add M$ and $\mathcal{P}(B\mbox{-}{\rm mod})$, where $B =
(\End_A M)^{op}$ ([ARS, p.33]).  If $M$ is {\it a generator} (i.e.,
$_AA\in {\rm add}M$), then we have

\vskip5pt

\begin{lem} \label{yoneda} \ Let $M$ be a generator of $A$-mod.
Then $\Hom_{A}(M, -): A\mbox{-}{\rm mod} \longrightarrow
B\mbox{-}{\rm mod}$ is fully faithful.
\end{lem}

\noindent {\bf Proof.}  Since $M$ is a generator, for any $X\in
A$-mod there is a surjective $A$-map $M^m \twoheadrightarrow X$ for
some positive integer $m$. This implies that $\Hom_{A}(M, -)$ is
faithful.

Let $X, Y \in A$-mod, and $f: \Hom_A(M, X) \longrightarrow \Hom_A(M,
Y)$ be  a
  $B$-map. By taking right $\add M$-approximations,
we get exact sequences $T_1 \xlongrightarrow{u} T_0
\xlongrightarrow{\pi} X \longrightarrow 0$ and $T_1^\prime
\xlongrightarrow{u^\prime} T_0^\prime \xlongrightarrow{\pi^\prime} Y
\longrightarrow 0$ with $T_0, T_1, T_0', T_1'\in\add M$ (since $M$
is a generator, $\pi$ and $\pi'$ are surjective). Applying
$\Hom_A(M, -)$ we have the following diagram with exact rows

  \[ \xymatrix{
  \Hom_{A}(M, T_1) \ar[r]^{\Hom_{A}(M,u)}  \ar@{-->}[d]^{f_1} &\Hom_{A}(M, T_0) \ar[r]^{\Hom_{A}(M,\pi)} \ar@{-->}[d]^{f_0} &\Hom_{A}(M, X)
  \ar[r] \ar[d]^f &0\\
\Hom_{A}(M, T_1^\prime) \ar[r]^{\Hom_{A}(M,u^\prime)} &\Hom_{A}(M,
T_0^\prime)
\ar[r]^{\Hom_{A}(M,\pi^\prime)}  &\Hom_{A}(M, Y) \ar[r] &0.\\
  }\]
Then  $f$ induces $f_1$ and $f_0$ such that the above diagram
commutes.
  Thus  $f_i = \Hom_{A}(M, f_i^\prime)$ for some $f_i^\prime\in\Hom_A(T_i, T_i^\prime)$,  $i = 0,
  1$.  So we get the following  diagram
  \[ \xymatrix{
  T_1 \ar[r]^u  \ar[d]^{f_1^\prime} &T_0 \ar[r]^\pi  \ar[d]^{f_0^\prime}
  &X
  \ar[r] \ar@{-->}[d]^{f^\prime} &0\\
T_1^\prime \ar[r]^{u^\prime}   &T_0^\prime \ar[r]^{\pi^\prime} &Y
  \ar[r]  &0\\
  }\]
with commutative left square.    So there exists
$f^\prime\in\Hom_A(X, Y)$ such that the above diagram commutes. Thus
$f\Hom_{A}(M,
  \pi) = \Hom_{A}(M, f^\prime)\Hom_{A}(M, \pi)$.  Since $\Hom_{A}(M, \pi)$ is surjective, it follows that $f = \Hom_{A}(M,
f^\prime)$, i.e., $\Hom_{A}(M, -)$ is full. $\s$

\vskip10pt

\subsection{} The following Auslander-Bridger
Lemma is very useful.

\begin{lem} \label{ABlemma} {\rm ([AB, Lemma
3.12])} \ Let $\mathcal A$ be an abelian category with enough
projective objects, $\mathcal X$ a resolving subcategory of
$\mathcal A$. Assume that
$$0\longrightarrow X_n\longrightarrow X_{n-1}\longrightarrow \cdots
\longrightarrow X_0 \longrightarrow A\longrightarrow 0$$ and
$$0\longrightarrow Y_n\longrightarrow Y_{n-1}\longrightarrow \cdots
\longrightarrow Y_0 \longrightarrow A\longrightarrow 0$$ are exact
sequences in $\mathcal A$, such that $X_i\in\mathcal X$ and
$Y_i\in\mathcal X$ for $0\le i\le n-1.$ Then $X_n\in\mathcal X$ if
and only if $Y_n\in\mathcal X$.\end{lem}

\noindent{\bf Proof.} \ We include a shorter proof . Assume
$X_{n}\in \mathcal{X}$. Take an exact sequence in $\mathcal A$
$$0\longrightarrow K\stackrel{d_{n}}\longrightarrow P_{n-1}\stackrel{d_{n-1}}\longrightarrow
 \cdots \longrightarrow P_{1}\stackrel{d_{1}}\longrightarrow P_{0}\stackrel{d_{0}}\longrightarrow A\longrightarrow 0$$
such that $P_{i}$ are projective, $i=0,\cdots,n-1.$ Then we get
chain map $f^{\bullet}$:
$$\CD P^{\bullet}: \quad@. 0@>>> K @>d_{n}>>  P_{n-1} @>d_{n-1}>>
\cdots @>>>P_{0}@>d_{0}>>A@>>>  0\cr @Vf^{\bullet}VV@. @V_{f_{n}}VV
@Vf_{n-1}VV   @. @Vf_{0}VV @| \cr X^{\bullet}:\quad@.0 @> >> X_{n}
@>\partial_{n}>> X_{n-1}@>\partial_{n-1}>> \cdots @> >>
  X_{0}@>\partial_{0}>>A@>>>0.
\endCD
$$

\vskip10pt

\noindent Consider the short exact sequence of complexes over
$\mathcal A$
$$0\longrightarrow X^\bullet \stackrel {\binom{0}{1}}\longrightarrow
{\rm Con}(f^\bullet) \stackrel{(1, 0)}\longrightarrow
P^\bullet[1]\longrightarrow 0,$$ where ${\rm Con}(f^\bullet)$
denotes the mapping cone of $f^\bullet$, which is by definition the
complex
$${\rm
Con}(f^\bullet): \quad 0 \longrightarrow K \stackrel{\begin{pmatrix}
  {-d}_{n}  \\
  f_{n} \\
\end{pmatrix}}{\longrightarrow} P_{n-1}\oplus
X_{n}\stackrel{\begin{pmatrix}
  {-d}_{n-1} & 0 \\
  f_{n-1}& \partial_{n} \\
\end{pmatrix}} \longrightarrow
P_{n-2}\oplus X_{n-1}$$

$$\longrightarrow \cdots \longrightarrow P_{0}\oplus X_{1}\stackrel{\begin{pmatrix}
  {-d}_{0} & 0 \\
  f_{0} & \partial_{1} \\
\end{pmatrix}}\longrightarrow
A\oplus X_{0} \stackrel{
  ({\rm Id}_{A}, \partial_{0} )}\longrightarrow A\longrightarrow 0.$$
Then ${\rm Con}(f^{\bullet})$ is again acyclic. Since $({\rm
Id}_{A},
\partial_{0})$ is a splitting epimorphism, it follows that we get an acyclic complex
$$ 0 \longrightarrow K
\longrightarrow P_{n-1}\oplus X_{n} \longrightarrow P_{n-2}\oplus
X_{n-1} \longrightarrow \cdots \longrightarrow P_{0}\oplus
X_{1}\stackrel {(f_0, \partial_1)}\longrightarrow
 X_{0} \longrightarrow 0$$
where $X_{0}\in\mathcal{X}$ and $P_{i-1}\oplus X_{i}\in\mathcal{X}$
for $i=1,\cdots,n$. Since $\mathcal{X}$ is closed under taking
kernels of epimorphisms, it follows that $K\in \mathcal{X}$.

\vskip10pt

By the similar way we get an acyclic complex
$$ 0 \longrightarrow K
\longrightarrow P_{n-1}\oplus Y_{n} \stackrel{\alpha}\longrightarrow
P_{n-2}\oplus Y_{n-1} \longrightarrow \cdots P_{0}\oplus
Y_{1}\longrightarrow
 Y_{0} \longrightarrow 0$$
where $Y_{0}\in\mathcal{X}$,   $P_{i-1}\oplus Y_{i}\in\mathcal{X}$
for $i=1, \cdots, n-1$, and ${\rm Im} \alpha\in \mathcal{X}$. Since
$\mathcal{X}$ is closed under taking extensions, by the short exact
sequence $0 \longrightarrow K \longrightarrow P_{n-1}\oplus Y_{n}
\longrightarrow{\rm Im }\alpha\longrightarrow 0$ with
$K\in\mathcal{X}$ and ${\rm Im} \alpha\in \mathcal{X}$, we know
$P_{n-1}\oplus Y_{n}\in \mathcal{X}$. Since $\mathcal{X}$ is closed
under taking direct summands, we have $Y_{n}\in \mathcal{X}$. $\s$

\vskip10pt

\subsection{} The following result is one of the key steps in our
categorical resolutions.

\begin{thm}\label{global} \ Let $A$ be an Artin algebra,  $T$ and $M$ modules in
$A\mbox{-}{\rm mod}$ such that  $^\perp T = {\rm add} M.$ Put $B: =
(\End_A M)^{op}$. Then for each positive integer $r\ge 2$,  ${\rm
gl.dim} B\le r$ if and only if  ${\rm inj.dim} T\le r$.\end{thm}
\noindent{\bf Proof.} \ Assume that ${\rm gl.dim}B\le r$. Let $X\in
A$-mod. Consider a right ${\rm add}M$-approximation $f_0:
M_0\longrightarrow X$. Since $M$ is a generator, $f_0$ is
surjective. Again considering a right ${\rm add}M$-approximation
$M_1\longrightarrow {\rm Ker}f_0$ and continuing this process we get
an exact sequence in $A$-mod
$$M_{r-1} \stackrel{f_{r-1}}\longrightarrow M_{r-2} \longrightarrow \cdots
\longrightarrow M_0\stackrel {f_0}\longrightarrow  X\longrightarrow
0$$ with $M_i\in {\rm add}M, \ 0\le i\le r-1.$ Put $K: = {\rm
Ker}f_{r-1}$. By construction we get an exact sequence
$$0\rightarrow \Hom(M, K)\longrightarrow \Hom(M,
M_{r-1})\longrightarrow \cdots \longrightarrow \Hom(M, M_0)
\longrightarrow \Hom(M, X) \rightarrow 0.$$ Since by assumption
${\rm proj.dim} _B\Hom(M, X)\le r$, it follows from
Auslander-Bridger Lemma that $\Hom_A(M, K)$ is projective as a
$B$-module. Thus there is a $B$-isomorphism $s: {\rm Hom}_A(M,
K)\longrightarrow {\rm Hom}_A(M, M')$ with $M'\in {\rm add}M$. By
Lemma \ref{yoneda} there exists $f: K\longrightarrow M'$ and $g:
M'\longrightarrow K$ such that $s = {\rm Hom}_A(M, f)$ and $s^{-1} =
{\rm Hom}_A(M, g).$ Therefore
$${\rm Hom}_A(M, fg) = {\rm Id}_{{\rm Hom}_A(M, M')} = {\rm Hom}_A(M, {\rm
Id}_{M'}).
$$
Since $\Hom_A(M, -)$ is faithful, it follows that $fg = {\rm
Id}_{M'}.$ Thus $K\in{\rm add}M$, and hence we have an exact
sequence in $A$-mod
$$0\longrightarrow M_r \longrightarrow M_{r-1} \longrightarrow \cdots
\longrightarrow M_0\longrightarrow  X\longrightarrow 0$$ with
$M_i\in {\rm add}M, \ 0\le i\le r.$ For $i\ge 1$, by dimension shift
we have
$$\Ext^{i+r}_A(X, T) \cong \Ext_A^{i}(M_r, T) =
0$$ since $M\in \ ^\perp T$. This shows  ${\rm inj.dim} T \le r$.

\vskip5pt

Note that the argument above holds also for $r\le 1$.

\vskip5pt

Conversely, assume that ${\rm inj.dim} T \le r$ with $r\ge 2$. Let
$_BY\in B$-mod. Taking a projective presentation of $_BY$
$$\Hom_A(M, M_1)\stackrel d \longrightarrow \Hom_A(M, M_0)
\longrightarrow \ _BY \longrightarrow 0$$ with $M_i\in {\rm add} M,
\ i = 0, 1$. Then there exists an $A$-map $f: M_1\longrightarrow
M_0$ such that $d = \Hom_A(M, f)$. Considering a right ${\rm
add}M$-approximation $M_2\longrightarrow {\rm Ker}f$ and continuing
this process, we get an exact sequence in $A$-mod
$$0\longrightarrow K \longrightarrow M_{r-1} \longrightarrow \cdots
\longrightarrow M_2\longrightarrow M_1 \stackrel f\longrightarrow
M_0$$ with $M_i\in {\rm add}M, \ 0\le i\le n-1.$ For $i\ge 1$, by
dimension shift we have
$$\Ext^i_A(K, T) \cong \Ext_A^{i+r-2}({\rm Ker}f, T) \cong
\Ext_A^{i+r-1}({\rm Im}f, T) \cong \Ext_A^{i+r}({\rm Coker}f, T) =
0$$ since ${\rm inj.dim} T \le r$. Thus $K\in \ ^\perp T = {\rm add}
M$. By construction we get an exact sequence \begin{align*}
0&\longrightarrow \Hom_A(M, K)\longrightarrow \Hom_A(M,
M_{r-1})\longrightarrow \cdots \longrightarrow \Hom_A(M, M_2) \\ &
\longrightarrow \Hom_A(M, M_1) \stackrel d\longrightarrow \Hom_A(M,
M_0) \longrightarrow \ _BY\longrightarrow 0.\end{align*} This gives
a projective resolution of $_BY$, and hence ${\rm proj.dim} _BY\le
r$. This proves ${\rm gl.dim}B\le r$. $\s$

\vskip10pt

\begin{rem} \ $(i)$ If $A$ is representation-finite and $T$ an injective module, then Theorem \ref{global}
is Auslander's result, and the corresponding $B$ is the Auslander
algebra {\rm([ARS])}.

\vskip5pt
 $(ii)$ If $A$ is {\rm CM}-finite Gorenstein algebra and  $T =\ _AA$, then Theorem \ref{global} is also well-known {\rm (}see
{\rm[LZ]} and {\rm [Bel2];} also {\rm [Leu])},  and the
corresponding $B$ is the relative Auslander algebra of $A$.
\end{rem}

\section{\bf Relative derived categories}

For our purpose we recall some points of relative derived
categories. Let $\mathcal C$ be a full additive subcategory of an
abelian category $\mathcal A$.

\vskip5pt

\subsection{} A complex $M^\bullet$ over $\mathcal A$ is {\it $\mathcal
C$-acyclic}, if $\operatorname{Hom}_{\mathcal A}(C, M^\bullet)$ is
acyclic for all objects $C$ in $\mathcal C$. A chain map \
$f^\bullet: X^\bullet \longrightarrow Y^\bullet$ is {\it a $\mathcal
C$-quasi-isomorphism}, if \ $\operatorname{Hom}_{\mathcal A}(C,
f^\bullet): \operatorname{Hom}_{\mathcal A}(C,
X^\bullet)\longrightarrow \operatorname{Hom}_{\mathcal A}(C,
Y^\bullet)$ is a quasi-isomorphism for all objects $C$ in $\mathcal
C$. Then $f^\bullet$ is a $\mathcal C$-quasi-isomorphism if and only
if the mapping cone \ $\operatorname{Con}(f^\bullet)$ \ is $\mathcal
C$-acyclic.

\vskip5pt

For $*\in \{b, -, {\rm blank}\}$, let  $K^*_{\mathcal C ac}(\mathcal
A)$ denote the full subcategory of the homotopy category
$K^*(\mathcal A)$ consisting of $\mathcal C$-acyclic complexes. Then
$$K^*_{\mathcal C ac}(\mathcal A) = \ ^\perp\mathcal C: = \{X^\bullet\in K^*(\mathcal A)\ | \
\Hom_{K^*(\mathcal A)}(C, X^\bullet[n]) =0, \ \forall \ n\in\Bbb Z,
\ \forall \ C\in\mathcal C\}.$$ Thus $K^*_{\mathcal C ac}(\mathcal
A)$ is a thick subcategory of $K^*(\mathcal A)$. The Verdier
quotient  $$D^*_{\mathcal C}(\mathcal A): = K^*(\mathcal
A)/K^*_{\mathcal Cac}(\mathcal A)$$ is called {\it the $\mathcal
C$-relative derived category}. See [GZ], [C2] and [AHV].

\vskip5pt

\begin{exm} \label{exmrelder}

$(i)$ \  If $\mathcal A$ has enough projective objects with
$\mathcal P$ the full subcategory consisting of projective objects,
and $\mathcal C = \mathcal {P}$, then $D^*_{\mathcal C}(\mathcal A)$
is just the derived category $D^*(\mathcal A)$.

\vskip5pt

$(ii)$ \ Let $\mathcal A$ be as in $(i)$. If $\mathcal C = \mathcal
{GP}(\mathcal A)$, the full subcategory of the Gorenstein-projective
objects of $\mathcal A$, then $D^*_{\mathcal C}(\mathcal A)$ is the
Gorenstein derived category in {\rm  [GZ]}.

\vskip5pt

$(iii)$ \ Let $A$ be an Artin algebra, and $M\in A$-mod. Then we
have the $M$-relative derived categories $D^*_{{\rm add}
M}(A\mbox{-}{\rm mod})$ and $D^*_{{\rm Add} M}(A\mbox{-}{\rm Mod})$.
See {\rm [AHV]}.
\end{exm}

\subsection{} It is important that the upper bounded derived category $D^-(\mathcal A)$ is a
triangulated subcategory of the unbounded derived category
$D(\mathcal{A})$, and that the bounded derived category
$D^{b}(\mathcal{A})$ is a triangulated subcategory of the
$D^{-}(\mathcal{A})$. The $\mathcal C$-relative derived category
enjoy this property. The proof is similar as the Gorenstein derived
category ([GZ, 2.5]), with a minor change. For the convenience of
the reader we include the proof.

\vskip5pt

\begin{lem} \ \label{conditionforsub} {\rm ([K2], Lemma
10.3)} \ Let $\mathcal B$ and $\mathcal D$ be triangulated
subcategories of triangulated category $\mathcal C$. If one of the
following conditions is satisfied, then the canonical triangle
functor $\mathcal D/\mathcal D\cap \mathcal B \longrightarrow
\mathcal C/\mathcal B$ is fully faithful.

\vskip5pt

$(i)$ \ Each morphism $X\longrightarrow B$ with $B \in\mathcal B$
and $X\in \mathcal D$ admits a factorization $X\longrightarrow
B'\longrightarrow B$ with $B'\in \mathcal D\cap \mathcal B$.

\vskip5pt

$(ii)$\ Each morphism $B\longrightarrow Y$ with $B \in\mathcal B$
and $Y\in \mathcal D$ admits a factorization $B\longrightarrow
B'\longrightarrow Y$ with $B'\in \mathcal D\cap \mathcal B$.
\end{lem}

\begin{prop} \label{reltivesub} \ Let $\mathcal C$ be a full additive subcategory of an
abelian category $\mathcal A$. Then $D^{-}_{\mathcal
C}(\mathcal{A})$ is a triangulated subcategory of $D_{\mathcal
C}(\mathcal{A})${\rm;} and $D^{b}_{\mathcal C}(\mathcal{A})$ is a
triangulated subcategory of $D^{-}_{\mathcal C}(\mathcal{A})$, and
hence of $D_{\mathcal C}(\mathcal{A})$.
\end{prop}
\noindent {\bf Proof.} \ We prove the first assertion by Lemma
\ref{conditionforsub}$(i)$, the second one can be proved by Lemma
\ref{conditionforsub}$(ii)$. Let $f^\bullet: \
X^\bullet\longrightarrow B^\bullet$ be a chain map with
$B^\bullet\in K_{\mathcal Cac}(\mathcal A)$ and $X^\bullet\in
K^-(\mathcal{A})$. We may assume that $X^i=0$ for $i> 0$. Then
$f^\bullet$ admits the following natural factorization:
\[\xymatrix{X^\bullet\ar[d]_{f^\bullet}: & \cdots  \ar [r] & X^{-1}
\ar[d] \ar[r] & X^0 \ar[d] \ar[r] & 0 \ar[r]\ar[d] & 0 \ar[d]
\ar[r]& \cdots \\
B'^\bullet:\ar[d] & \cdots  \ar [r] & B^{-1}\ar@{=}[d] \ar[r] & B^0
\ar@{=}[d] \ar[r] &
\operatorname{Ker}d^1 \ar[r]\ar[d] & 0 \ar[r]\ar[d] &\cdots \\
B^\bullet: & \cdots  \ar[r] & B^{-1} \ar[r] & B^0\ar[r] & B^1 \ar[r]
& B^2 \ar[r] & \cdots }\]

\noindent We need to prove that $B'^\bullet$ is $\mathcal
C$-acyclic. Since $B^\bullet$ is $\mathcal C$-acyclic, it suffices
to prove that
$$\operatorname{Hom}_\mathcal A(C,
B^{-1})\stackrel{d^{-1}_*}{\longrightarrow}
\operatorname{Hom}_\mathcal A(C, B^0) \stackrel {\widetilde{d^0}_*}
 {\longrightarrow} \operatorname{Hom}_\mathcal A(C,
\operatorname{Ker}d^1)\longrightarrow 0$$ is exact for each $C\in
\mathcal C$, where $\widetilde{d^0}: B^0\longrightarrow
\operatorname{Ker}d^1$ is induced by $d^0$.  Since
$0\longrightarrow\operatorname{Hom}(C,
\operatorname{Ker}d^1)\stackrel \sigma\hookrightarrow
\operatorname{Hom}(C, B^1)\longrightarrow \operatorname{Hom}(C,
B^2)$ is exact, by the commutative diagram
\[\xymatrix{\operatorname{Hom}(C, B^{-1})
\ar[r]^-{d^{-1}_*}& \operatorname{Hom}(C, B^0)
 \ar[rd]_- {\widetilde {d^0}_*}
 \ar[rr]^-{d^0_*}
&&\operatorname{Hom}(C, B^1)\ar[r]^-{d^1_*}& \operatorname{Hom}(C, B^2)\\
&& \operatorname{Hom}(C, \operatorname{Ker}d^1)\ar
@^{(->}[ru]_\sigma }\]  we have $\operatorname{Ker}\widetilde{d^0}_*
= \operatorname{Ker}d_0^* = \operatorname{Im}d^{-1}_*$, and
$\operatorname{Im}\widetilde {d^0}_* = \operatorname{Im} d^0_* =
\operatorname{Ker}d^1_* = \operatorname{Hom}_\mathcal A(C,
\operatorname{Ker}d^1).$ \hfill$\blacksquare$

\vskip10pt

The natural functor $\mathcal{A} \longrightarrow D^b_{\mathcal
C}(\mathcal{A})$, which is the composition of the embedding
$\mathcal{A} \longrightarrow K^b(\mathcal A)$ and the localization
functor $K^b(\mathcal A)\longrightarrow D^b_{\mathcal
C}(\mathcal{A})$, is fully faithful. The proof is as [GZ, 2.9].

\vskip10pt

\subsection{} If $\mathcal A$ has enough projective objects and $\mathcal
P\subseteq \mathcal C$, then $K^*_{\mathcal Cac}(\mathcal A)$ is a
thick subcategory of $K^*_{ac}(\mathcal A)$, where
$K^*_{ac}(\mathcal A)$ is the full subcategory of the homotopy
category $K^*(\mathcal A)$ consisting of acyclic complexes. By Lemma
\ref{Verdier} there is a triangle-equivalence
$$D^*(\mathcal A)\cong D^*_{\mathcal C}(\mathcal A)
/(K^*_{ac}(\mathcal A)/K^*_{\mathcal Cac}(\mathcal A)),$$ and we
have the localization functor $\pi_*: D^*_{\mathcal C}(\mathcal
A)\longrightarrow D^*(\mathcal A)$. Note that $\pi_*$ is an
equivalence if and only if $\mathcal {C} = \mathcal P$.

\vskip10pt

\begin{lem} \label{basic} \  Let $\mathcal C$ be a full additive subcategory of abelian category $\mathcal A$.  Then we have

\vskip5pt

$(i)$ \ {\rm ([CFH], Proposition 2.6)} \ A chain map $f^\bullet:
X^\bullet \longrightarrow Y^\bullet$ \ is a
$\mathcal{C}$-quasi-isomorphism if and only if there are
isomorphisms of abelian groups for any $C^\bullet\in K^-(\mathcal
C):$
$$\operatorname{Hom}_{K(\mathcal A)}(C^\bullet, f^\bullet[n]):
\ \operatorname{Hom}_{K(\mathcal A)}(C^\bullet, X^\bullet[n])\cong
\operatorname{Hom}_{K(\mathcal A)}(C^\bullet, Y^\bullet[n]), \
\forall \ n\in\Bbb Z.$$

\vskip5pt

$(ii)$  \ {\rm ([GZ], Lemma 2.2)} \  Let $C^\bullet\in K^-(\mathcal
C)$, and \ $f^\bullet:  X^\bullet \longrightarrow C^\bullet$ be a
$\mathcal{C}$-quasi-isomorphism. Then there is $g^\bullet: C^\bullet
\longrightarrow X^\bullet$ such that $f^\bullet g^\bullet$ is
homotopic to ${\rm Id}_{C^\bullet}$.

\vskip5pt

Thus, if in addition $X^\bullet\in K^-(\mathcal C)$, then
$f^\bullet$ is a homotopy equivalence.

\vskip5pt

$(iii)$  \ {\rm ([GZ], Proposition 2.8)} \  Let $C^\bullet\in
K^-(\mathcal C)$ and $Y^\bullet$ be an arbitrary complex. Then $Q: \
f^\bullet\mapsto f^\bullet/{\rm Id}_{C^\bullet}$ \ gives an
isomorphism
 \ $\operatorname{Hom}_{K(\mathcal{A})}(C^\bullet, Y^\bullet)\cong
\operatorname{Hom}_{D_{\mathcal C}(\mathcal{A})}(C^\bullet,
Y^\bullet)$ of abelian groups.

\vskip5pt

In particular, $K^b(\mathcal C)$ can be viewed as a triangulated
subcategory of $D^b_{\mathcal C}(\mathcal{A});$ and $K^-(\mathcal
C)$ can be viewed as a triangulated subcategory of $D^-_{\mathcal
C}(\mathcal{A})$.
\end{lem}

\subsection{} Let $K^{-, \mathcal C b}(\mathcal C)$ denote the full
subcategory of $K^{-}(\mathcal C)$ given by
\begin{align*}K^{-, \mathcal C b}(\mathcal C):=
\{X^\bullet\in K^{-}(\mathcal C)\ | & \ \exists \ N\in\Bbb Z \
\mbox{such that}\ {\rm H}^i\operatorname{Hom}_\mathcal A(C,
X^\bullet) = 0, \\ & \ \ \forall \ i\le N, \ \forall \ C\in\mathcal
C\}.\end{align*} Then  $K^{-, \mathcal C b}(\mathcal C)$ is a thick
triangulated subcategory of $K^{-}(\mathcal C)$.

\vskip10pt

\begin{lem} \label{gpqis} \ Let $\mathcal C$ be a contravariantly finite subcategory of abelian category $\mathcal A$.
Then for each $X^\bullet \in K^b(\mathcal A)$ there is a
$\mathcal{C}$-quasi-isomorphism $C_{_{X^\bullet}} \longrightarrow
X^\bullet$  with $C_{_{X^\bullet}}\in K^{-, \mathcal C b}(\mathcal
C)$.
\end{lem}
\noindent{\bf Proof.} \ The proof is similar as [GZ, Proposition
3.4]. Use induction on the width $w(X^{\bullet})$, the number of $i$
such that $X^i\ne 0$. Assume that $w(X^{\bullet}) = 1$. Then
$X^{\bullet}$ is the stalk complex of object $X$, say at degree $0$.
Since $\mathcal C$ is contravariantly finite in $\mathcal A$, there
exists a right $\mathcal {C}$-approximation $d^0: C^0
\longrightarrow X$ of $X$. Taking a right $\mathcal
{C}$-approximation $C^{-1} \longrightarrow {\rm Ker}d^0$, and
continuing this process we get a complex
$$C^\bullet: \cdots \longrightarrow C^{-1} \stackrel
{d^{-1}}\longrightarrow C^0\stackrel {d^{0}}\longrightarrow
X\longrightarrow 0$$ with each $C^{i}\in \mathcal {C}$, such that
$\Hom_\mathcal A(C, C^\bullet)$ is acyclic for each $C\in \mathcal
{C}$. Put $C_{_{X^\bullet}}$ to be the complex obtained from
$C^\bullet$ by deleting $X$. By construction we get a
$\mathcal{C}$-quasi-isomorphism $\phi_{_{X^\bullet}}: \
C_{_{X^\bullet}} \longrightarrow X^\bullet$ with
$C_{_{X^\bullet}}\in K^{-, \mathcal C b}(\mathcal C)$.

\vskip5pt

Assume $w(X^{\bullet})\ge 2$ with $X^j\ne 0$ and $X^i = 0$ for $i <
j$. Then we have a distinguished triangle \
$X^\bullet_1\stackrel{u}{\longrightarrow} X^\bullet_2\longrightarrow
X^\bullet\longrightarrow X^\bullet_1[1]$ in $K^b(\mathcal A)$, where
$X^\bullet_1:= X^j[-j-1]$ and $X^\bullet_2$ is the brutal truncated
complex $X^{\bullet}_{>j}.$ By induction there exist
$\mathcal{C}$-quasi-isomorphisms
$$\phi_1: C_{_{X^\bullet_1}}\longrightarrow X^\bullet_1, \
\ \ \phi_2: C_{_{X^\bullet_2}}\longrightarrow X^\bullet_2$$ with
$C_{_{X_1^\bullet}}, \ C_{_{X_2^\bullet}}\in K^{-, \mathcal C
b}(\mathcal C)$. By Lemma \ref{basic}$(i)$ $\phi_2$ induces an
isomorphism
$$\operatorname{Hom}_{K^{-}(\mathcal A)}(C_{_{X^\bullet_1}}, \ C_{_{X^\bullet_2}})
\cong \operatorname{Hom}_{K^{-}(\mathcal A)}(C_{_{X^\bullet_1}}, \
X^\bullet_2).$$ Thus there is a unique chain map \ $f^\bullet:
C_{_{X^\bullet_1}}\longrightarrow C_{_{X^\bullet_2}}$ \ such that
$\phi_2\circ f^\bullet = u\circ \phi_1$. Embedding $f^\bullet$ into
a distinguished triangle in $K^{-, \mathcal C b}(\mathcal C)$
$$C_{_{X^\bullet_1}}\stackrel {f^\bullet} {\longrightarrow}C_{_{X^\bullet_2}}\longrightarrow
C_{_{X^\bullet}} \longrightarrow C_{_{X^\bullet_1}}[1]$$ we get a
unique complex $C_{_{X^\bullet}}$ in $K^{-, \mathcal C b}(\mathcal
C)$. By the axiom of a triangulated category, there is
$\phi_{_{X^\bullet}}: C_{_{X^\bullet}}\longrightarrow X^\bullet$
such that the diagram commutes
$$
\xymatrix{
   C_{_{X^\bullet_1}}\ar[r]^{f^\bullet}\ar[d]_{\phi_1} & C_{_{X^\bullet_2}} \ar[r]\ar[d]_{\phi_2} & C_{_{X^\bullet}}
\ar[r]\ar@{..>}[d]_{\phi_{_{X^\bullet}}}
  & C_{_{X^\bullet_1}}[1] \ar[d]_{\phi_1[1]} \\
   X^\bullet_1\ar[r]^{u} & X^\bullet_2\ar[r]& X^\bullet \ar[r] & X^\bullet_1[1].}$$
By using cohomological functors and the Five-Lemma it is easy to
know that $\phi_{_{X^\bullet}}$ \ is a
$\mathcal{C}$-quasi-isomorphism. $\s$

\vskip 10pt

The following result is due to J.Asadollahi, R.Hafezi, and R.Vahed
[AHV, Theorem 3.3] (see also [GZ, Theorem 3.6] for the Grorenstein
derived category). Since we
 need the equivalence $F: K^{-, \mathcal
Cb}(\mathcal C)\longrightarrow D^b_{\mathcal C}(\mathcal A)$ in its
proof, and since the proof was omitted in [AVH], so we include a
proof.

\vskip10pt

\begin{prop} \label{gpder} {\rm ([AHV])} \ Let $\mathcal C$ be a contravariantly finite subcategory of abelian category $\mathcal A$.
Then there is a triangle-equivalence \ $D^b_{\mathcal C}(\mathcal
A)\cong K^{-, \mathcal Cb}(\mathcal C),$ which fixes objects in
$K^b(\mathcal C)$.
\end{prop} \noindent{\bf Proof.} \ Let $F: K^{-, \mathcal Cb}(\mathcal
C)\longrightarrow D^-_{\mathcal C}(\mathcal A)$ be the composite of
the embedding $K^{-, \mathcal Cb}(\mathcal C)\hookrightarrow
K^{-}(\mathcal A)$ and the localization functor $Q: K^{-}(\mathcal
A)\longrightarrow D^-_{\mathcal C}(\mathcal A)$. For each complex
$X^\bullet\in K^{-, \mathcal Cb}(\mathcal C)$, by definition there
is an $N\in\Bbb Z$ such that ${\rm H}^{i}\operatorname{Hom}_\mathcal
A(C, X^\bullet) = 0, \ \forall \ i \le N, \ \forall \
C\in\mathcal{C}$. Since the following chain map is a
$\mathcal{C}$-quasi-isomorphism
\[\xymatrix{X^\bullet\ar[d]_{f^\bullet}:
& \cdots \ar [r] & X^{N-2} \ar [r] \ar[d]& X^{N-1} \ar[d] \ar[r] &
X^{N} \ar@{=}[d] \ar[r] & X^{N+1}\ar[r]\ar@{=}[d]& \cdots
 \\
\tau_{\ge N}X^\bullet: & \cdots  \ar [r] & 0 \ar [r] &
\operatorname{Ker}d^{N} \ar[r] & X^N \ar[r] & X^{N+1}\ar[r]&
 \cdots}\] it follows that there is an isomorphism $F(X^\bullet)\cong
\tau_{\ge N}X^\bullet$ in $D^-_{\mathcal C}(\mathcal A)$ with
$\tau_{\ge N}X^\bullet\in D^b_{\mathcal C}(\mathcal A)$. Thus the
image of $F$ falls in $D^b_{\mathcal C}(\mathcal A)$, and hence $F$
induces a triangle functor $K^{-, \mathcal Cb}(\mathcal
C)\longrightarrow D^b_{\mathcal C}(\mathcal A)$, again denoted by
$F$ (here we need to use Proposition \ref{reltivesub}). In
particular, $F$ fixes objects in  $K^b(\mathcal C)$, i.e.,
$F(X^\bullet) = X^\bullet$ for $X^\bullet\in  K^b(\mathcal C)$.

\vskip5pt

By Lemma \ref{gpqis} $F$ is dense; and by Lemma \ref{basic}$(iii)$
$F$ is fully faithful. $\s$

\vskip10pt

\section{\bf A relative description of bounded derived category}

\subsection {} Let $\mathcal A$ be an abelian category with enough projective objects, and $\mathcal C$ a resolving subcategory of $\mathcal A$.
By $K^b_{ac}(\mathcal C)$ we denote the full subcategory of
$K^-(\mathcal C)$ consisting of those complexes which are homotopy
equivalent to bounded acyclic complexes over $\mathcal C$. It is
clear that $K^b_{ac}(\mathcal C)$ is a triangulated subcategory of
$K^-(\mathcal C)$.

\begin{lem} \label{b}
Let $C^\bullet\in K^{-, \mathcal C b}(\mathcal C)$. If $C^\bullet$
is acyclic, then $C^\bullet \in K^b_{ac}(\mathcal C)$.
\end{lem} \noindent {Proof.} \ Since $C^\bullet = (C^i,
d^i)$ is upper bonded acyclic complex over $\mathcal C$, and
$\mathcal C$ is closed under kernels of epimorphisms, it follows
that $\Ima{d^i} \in \mathcal C, \ \forall \ i\in \Bbb Z$. Since
$C^\bullet\in K^{-, \mathcal C b}(\mathcal C)$, by definition there
exists an integer $N$ such that ${\rm H}^{n}\Hom_{\mathcal A}(C,
C^\bullet) = 0, \ \forall \ n \le N, \ \forall \ C \in \mathcal C$.
In particular ${\rm H}^{n}\Hom_{\mathcal A}(\Ima{d^{n-1}},
C^\bullet) = 0.$ This implies that the induced epimorphism
$\widetilde{d^{n-1}}: C^{n-1} \longrightarrow \Ima d^{n-1}$ splits
for $n \le N$, and hence there is an isomorphism $C^\bullet\cong
C^{'\bullet}$ in $K^-(\mathcal C)$, where $C'^\bullet$ is the
complex $$\cdots \longrightarrow 0 \longrightarrow \Ima d^{N-1}
\hookrightarrow C^{N} \longrightarrow C^{N+1}\longrightarrow
\cdots$$ with  $C'^\bullet\in K^b_{ac}(\mathcal C).$ Thus $C^\bullet
\in K^b_{ac}(\mathcal C)$. $\s$

\subsection {} The following result is another key step in proving Theorem
\ref{mainresult}, and also it seems to be of independent interest.
If one takes $\mathcal C$ to be $\mathcal P$, then it read as the
well-known triangle equivalence $D^b(\mathcal A)\cong K^{-,
b}(\mathcal P)$. If $\mathcal C= \mathcal {GP}(\mathcal A),$ then it
is Theorem 5.1 of [KZ].

\vskip10pt

\begin{thm}\label{relativedescription} Let $\mathcal A$ be an abelian category with enough projective objects, and $\mathcal C$ a resolving contravariantly
finite subcategory of $\mathcal A$. Then $K^b_{ac}(\mathcal C)$ is a
thick subcategory of $K^{-,\mathcal C b}(\mathcal C)$, and we have a
triangle-equivalence
$$G: D^b(\mathcal A)\longrightarrow K^{-,\mathcal C b}(\mathcal C)/K^b_{ac}(\mathcal C)$$
such that $G$ sends an object $C\in K^b(\mathcal C)$ to $C\in
K^{-,\mathcal C b}(\mathcal C)/K^b_{ac}(\mathcal C).$
\end{thm}
\noindent{\bf Proof.} Lemma \ref{b} implies that $K^b_{ac}(\mathcal C)$ is a
thick subcategory of $K^{-,\mathcal C b}(\mathcal C)$.

\vskip5pt

Let $F': K^b_{ac}(\mathcal C)\longrightarrow K^b_{ac}(\mathcal
A)/K^b_{\mathcal Cac}(\mathcal A)$ be the composite of the embedding
functor $K^b_{ac}(\mathcal C)\hookrightarrow K^{b}_{ac}(\mathcal A)$
and the Verdier functor $Q: K^{b}_{ac}(\mathcal A)\longrightarrow
K^b_{ac}(\mathcal A)/K^b_{\mathcal Cac}(\mathcal A)$. We first claim
that $F'$ is a triangle equivalence.

\vskip5pt

Since $K^b_{ac}(\mathcal A)$ is a triangulated subcategory of
$K^b(\mathcal A)$, it follows that $K^b_{ac}(\mathcal
A)/K^b_{\mathcal Cac}(\mathcal A)$ is a triangulated subcategory of
$K^b(\mathcal A)/K^b_{\mathcal Cac}(\mathcal A)$. By definition
$K^b(\mathcal A)/K^b_{\mathcal Cac}(\mathcal A)$ is the $\mathcal
C$-relative derived category $D^b_\mathcal C(\mathcal A)$.  By Lemma
\ref{basic}$(iii)$ $F'$ is fully faithful.

\vskip5pt

For each complex $X^\bullet\in K^{b}_{ac}(\mathcal A)$, by Lemma
\ref{gpqis} there is a $\mathcal C$-quasi-isomorphism
$C^\bullet\longrightarrow X^\bullet$  with $C^\bullet\in K^{-,
\mathcal Cb}(\mathcal C)$. Since $\mathcal C\supseteq \mathcal P$,
it follows that a $\mathcal C$-quasi-isomorphism is a
quasi-isomorphism. Since $X^\bullet$ is acyclic, it follows that
$C^\bullet$ is acyclic. By Lemma \ref{b} $C^\bullet \in
K^b_{ac}(\mathcal C)$. By $X\cong F'(C^{\bullet})$ in
$K^b_{ac}(\mathcal A)/K^b_{\mathcal Cac}(\mathcal A)$ with
$C^\bullet \in K^b_{ac}(\mathcal C)$ we know that $F'$ is dense.
This proves the claim.

\vskip5pt

By construction $F'$ is just the restriction of $F$ to
$K^b_{ac}(\mathcal C)$, where $F$ is the triangle-equivalence $K^{-,
\mathcal Cb}(\mathcal C) \longrightarrow D^b_{\mathcal C}(\mathcal
A): = K^b(\mathcal A)/K^b_{\mathcal C ac}(\mathcal A)$ given in the
proof of Proposition \ref{gpder}. Hence we have a commutative
diagram
\[\xymatrix {K^b_{ac}(\mathcal C)\ar[r]\ar[d] & K^{-,
\mathcal C b}(\mathcal C)\ar[d]
\\ K^b_{ac}(\mathcal
A)/K^b_{\mathcal C ac}(\mathcal A)\ar[r] & K^b(\mathcal
A)/K^b_{\mathcal C ac}(\mathcal A)}\] where the horizontal functors
are embeddings, and the vertical ones are triangle-equivalences.
Thus $F$ induces a triangle-equivalence
$$K^{-,\mathcal Cb}(\mathcal C)/K^b_{ac}(\mathcal C)\cong (K^b(\mathcal A)/K^b_{\mathcal Cac}(\mathcal A))/(K^b_{ac}(\mathcal
A)/K^b_{\mathcal Cac}(\mathcal A))\eqno(*)$$ While by Lemma
\ref{Verdier} we have a triangle-equivalence $$\mbox{the right hand
side of } (*) \cong K^b(\mathcal A)/K^b_{ac}(\mathcal A) =
D^b(\mathcal A).$$ This proves $D^b(\mathcal A)\cong K^{-,\mathcal C
b}(\mathcal C)/K^b_{ac}(\mathcal C).$ Since this equivalence is
induced by $F$, and $F$ fixes objects in  $K^b(\mathcal C)$ by
Proposition \ref{gpder}, it follows that it sends an object $C\in
K^b(\mathcal C)$ to $C\in K^{-,\mathcal C b}(\mathcal
C)/K^b_{ac}(\mathcal C).$ $\s$

\vskip10pt

From the proof above, we see that the assumption ``$\mathcal C$ is a
resolving contravariantly finite subcategory" is used.

\vskip10pt

\subsection{} We need the following fact in the next
section.

\vskip10pt

\begin{prop} \label{ff} \ Let $\mathcal A$ be an abelian
category with enough projective objects, and $\mathcal C$ a full
additive subcategory of $\mathcal A$ with $\mathcal P\subseteq
\mathcal C$. Then there is a functorial isomorphism of abelian
groups for each $P\in K^-(\mathcal P)$ and $C\in K^{-}(\mathcal C)$
$$\Hom_{K^{-}(\mathcal C)}(P, C)\cong \Hom_{K^{-}(\mathcal C)/K^{b}_{ac}(\mathcal C)}(P,
C)$$ given by $f \mapsto f/{\rm Id}_P, \ \forall \ f\in
\Hom_{K^{-}(\mathcal C)}(P, C).$\end{prop}

\noindent{\bf Proof.} The proof is similar as in the case of derived
category. Since this assertion will be used, for the completeness we
include a justification.

\vskip5pt

Recall a well-known fact: if $t: Z \longrightarrow P$ is a
quasi-isomorphism with $P\in K^-(\mathcal P)$, then there is $g: P
\longrightarrow Z$ such that $t g$ is homotopic to ${\rm Id}_P$ (cf.
Lemma \ref{basic}$(ii)$).

\vskip5pt

Now assume $f/{\rm Id}_P = 0$. By definition we have a commutative
diagram in $K^{-}(\mathcal C)$
\[\xymatrix@R=0.9pc{
& & & & & P \ar@{=>}[ld]_-= \ar[rd]^-f &\\
& & & & P & Z \ar@{=>}[l]_-t \ar[d]^-t \ar[u]_-t \ar[r]^-{0} & C  & & & & \\
& & & & & P \ar@{=>}[lu]^-= \ar[ru]_-0}\] where $t: Z\longrightarrow
P$ a chain map such that ${\rm Con}(t)\in K^b_{ac}(\mathcal C).$
Thus $t$ is a quasi-isomorphism, and hence there is $g: P
\longrightarrow Z$ such that $t g$ is homotopic to ${\rm Id}_{P}$.
Thus by $ft = 0$ we have $f = f(t g) = 0$.

\vskip5pt

Assume  $f/s\in \Hom_{K^-(\mathcal C)/K^b_{ac}(\mathcal C)}(P, C)$,
where $s: Z\longrightarrow P$ with ${\rm Con}(s)\in
K^b_{ac}(\mathcal C)$, and $f: Z\longrightarrow C$. Since $s$ is
quasi-isomorphism, there is $g: P \longrightarrow Z$ such that $s g$
is homotopic to ${\rm Id}_{P}$, and hence we get a commutative
diagram
\[\xymatrix@R=0.9pc{
& & & & & Z \ar@{=>}[ld]_-s \ar[rd]^-f &\\
& & & & P & P \ar@{=>}[l]_-= \ar[d]^-= \ar[u]_-g \ar[r]^-{fg} & C  & & & & \\
& & & & & P \ar@{=>}[lu]^-= \ar[ru]_-{fg}}\] This means $f/s =
fg/{\rm Id}_P$. $\s$

\vskip10pt

\subsection{} For later use,  we need to investigate $K^{-}(\mathcal
C)/K^{b}_{ac}(\mathcal C)$ in more details.

\vskip5pt

Let $\mathcal A$ be an abelian category with enough projective
objects, and  $\mathcal C$ a resolving subcategory of $\mathcal A$.
An object $I\in \mathcal C$ is {\it a {\rm(}relative{\rm)} injective
object} of $\mathcal C$, provided that the functor $\Hom_\mathcal
A(-, I)$ sends any short exact sequence $0\longrightarrow
X_1\longrightarrow X_2\longrightarrow X_3\longrightarrow 0$ with
$X_i\in\mathcal C, \ i=1,2,3,$ to an exact sequence. Clearly, $I$ is
an injective object of $\mathcal C$ if and only if $\Ext_\mathcal
A^1(X, I) = 0$ for each $X\in\mathcal C$, also if and only if
$\Ext_\mathcal A^i(X, I) = 0$ for each $X\in\mathcal C$ and for
$i\ge 1$,

\vskip5pt

\begin{lem} \label{inj1} \  Let $\mathcal C$ be a resolving subcategory of $\mathcal
A$, and $G = (G^i, d^i_G)\in K^-_{ac}(\mathcal C)$. Assume that $I =
(I^i, d^i_I)$ is a bounded complex such that all $I^i$ are injective
objects of $\mathcal C$. Then $\Hom_{K^-(\mathcal A)}(G, I) = 0.$
\end{lem}
\noindent{\bf Proof.} The proof is similar to the case of $\mathcal
C = \mathcal A$, which is well-known. For the completeness we
include a proof.

Let $f: G\longrightarrow I$ be a chain map. We need to show that $f$
is null-homotopic. We construct a homotopy $s = (s^i)$ by induction.
Assume that we have constructed $s^i: G^i \longrightarrow I^i$ for
$i\le m$, such that $f^{i-1} = d^{i-2}_Is^{i-1} + s^id_G^{i-1}$ for
$i\le m$. Since $G$ is an upper bounded acyclic complex with all
$G^i\in \mathcal C,$ and $\mathcal C$ is closed under the kernels of
epimorphisms, it follows that ${\rm Im} d^j_G\in \mathcal C, \
\forall \ j\in\Bbb Z.$ Since
$$(f^m-d_I^{m-1}s^m)d_G^{m-1} = 0$$
it follows that $f^m-d_I^{m-1}s^m$ factors through ${\rm
Coker}d_G^{m-1} = {\rm Im}d^m_G.$ Since $I^m$ is an injective object
of $\mathcal C$, it follows that there is $s^{m+1}:
G^{m+1}\longrightarrow I^m$ such that
$$f^m-d_I^{m-1}s^m = s^{m+1}d^m_G.$$
This completes the proof.  $\s$

\vskip10pt

\begin{lem} \label{inj2} \  Let $\mathcal C$ be a resolving subcategory of $\mathcal
A$, and $C\in K^-(\mathcal C)$.  Assume that $I$ is a bounded
complex such that all $I^i$ are injective objects of $\mathcal C$.
If $t: I \longrightarrow C$ a quasi-isomorphism, then there exists a
chain map $s: C\longrightarrow I$ such that $st = {\rm Id}_I$ in
$K^-(\mathcal A)$.
\end{lem}
\noindent{\bf Proof.}  By Lemma \ref{inj1} $\Hom_{K^-(\mathcal
A)}({\rm Con} (t), I) = 0 = \Hom_{K^-(\mathcal A)}({\rm Con}(t)[-1],
I).$ Applying $\Hom_{K^-(\mathcal A)}(-, I)$ to the distinguished
triangle $I \stackrel{t}\longrightarrow C\longrightarrow {\rm Con}
(t) \longrightarrow I[1]$ we see that $\Hom_{K^-(\mathcal A)}(C,
I)\stackrel {\Hom(t, I)}\longrightarrow \Hom_{K^-(\mathcal A)}(I,
I)$ is an isomorphism, from which the assertion follows.  $\s$

\vskip10pt

\begin{prop} \label{inj3} \ Let $\mathcal A$ be an abelian
category with enough projective objects, and $\mathcal C$ a
resolving subcategory of $\mathcal A$. Assume that $I$ is a bounded
complex such that all $I^i$ are injective objects of $\mathcal C$.
Then for each $C\in K^{-}(\mathcal C)$ there is a functorial
{\rm(}in $C$ and in $I${\rm)} isomorphism
$$\Hom_{K^{-}(\mathcal C)}(C, I)\cong \Hom_{K^{-}(\mathcal C)/K^{b}_{ac}(\mathcal C)}(C,
I).$$\end{prop} \noindent{\bf Proof.} Here we need to use the left
fraction construction of $K^{-}(\mathcal C)/K^{b}_{ac}(\mathcal C)$.
The isomorphism is given by $f \mapsto {\rm Id}_I\backslash f, \
\forall \ f\in \Hom_{K^{-}(\mathcal C)}(C, I).$ The proof is dual to
the one of Proposition \ref{ff}, by using Lemma \ref{inj2}. We omit
the details. $\s$

\vskip10pt

\section{\bf Main results}

\subsection{} Now we are in position to prove

\vskip10pt

\begin{thm} \label{mainresult}  Let $A$ be an Artin
algebra with ${\rm gl.dim}A = \infty$. Assume that there are modules
$T$ and $M$ in $A\mbox{-}{\rm mod}$ with ${\rm inj.dim} T< \infty$,
such that $^\perp T = {\rm add}M$. Then $D^b(A\mbox{-}{\rm mod})$
admits a categorical resolution $D^b(B\mbox{-}{\rm mod})$ with  $B =
(\End_AM)^{op}$.
\end{thm}
\noindent{\bf Proof.} By Theorem \ref{global} ${\rm gl.dim} B <
\infty,$ i.e.,  $D^b(B\mbox{-}{\rm mod})$ is smooth.

\vskip5pt

The equivalence $\Hom_A(M, -): {\rm add M} \longrightarrow \mathcal
P(B\mbox{-}{\rm mod})$ of categories induces pointwisely a
triangle-equivalence $K^{-, {\rm add M}\ b}({\rm add M})\cong K^{-,
b}(\mathcal P(B\mbox{-}{\rm mod}))$. Since $D^b(B\mbox{-}{\rm mod})
\cong K^{-, b}(\mathcal P(B\mbox{-}{\rm mod}))$, we have a
triangle-equivalence
$$F: D^b(B\mbox{-}{\rm mod}) \cong K^{-, {\rm add M}\ b}({\rm add M}).$$
Since ${\rm add M} = \ ^\perp T$, it follows that ${\rm add M}$ is a
resolving contravariantly finite subcategory of $A$-mod, and hence
by Theorem \ref{relativedescription} we have a triangle-equivalence
$$G: D^b(A\mbox{-}{\rm mod})\longrightarrow K^{-, {\rm add M}\ b}({\rm add M})/K^{b}_{ac}({\rm add M})$$
such that $G$ sends an object $P\in K^b(\mathcal P(A\mbox{-}{\rm
mod}))$ to $P\in K^{-, {\rm add M}\ b}({\rm add M})/K^{b}_{ac}({\rm
add M}),$ i.e., $GP = P$. Thus, we get a triangle functor
$$\pi_*:= G^{-1}V F:  D^b(B\mbox{-}{\rm mod})  \longrightarrow
D^b(A\mbox{-}{\rm mod}),$$ where $V: K^{-, {\rm add M}\ b}({\rm add
M})\longrightarrow K^{-, {\rm add M}\ b}({\rm add
M})/K^{b}_{ac}({\rm add M})$ is the Verdier functor.

\vskip5pt

On the other hand,  by Proposition \ref{finitely filtrated}
$D^b(A\mbox{-}{\rm mod})^{\rm perf} = K^b(\mathcal P(A\mbox{-}{\rm
mod}))$. Thus we have a triangle functor
$$\pi^*:= F^{-1}\sigma:  D^b(A\mbox{-}{\rm mod})^{\rm perf} \longrightarrow
D^b(B\mbox{-}{\rm mod})$$ where $\sigma$ is the embedding $\sigma:
K^b(\mathcal P(A\mbox{-}{\rm mod}))\hookrightarrow K^{-, {\rm add
M}\ b}({\rm add M})$.

\vskip5pt

The diagram \[\xymatrix {D^{b}(B\mbox{-}{\rm
mod})\ar[d]^-{F}\ar[r]^{\pi_*} & D^{b}(A\mbox{-}{\rm mod})\ar[d]^-G
\\K^{-, {\rm add M}\ b}({\rm add
M})\ar[r]^-{V} & K^{-, {\rm add M}\ b}({\rm add M})/K^{b}_{ac}({\rm
add M})}\] commutes. Since $K^b_{ac}({\rm add M})$ is thick in
$K^{-,{\rm add M} b}(\mathcal C)$ (cf. Theorem
\ref{relativedescription}), we have ${\rm Ker} V = K^b_{ac}({\rm add
M})$. It follows that ${\rm Ker} \pi_* = F^{-} (K^{b}_{ac}({\rm add
M}))$, and $\pi_*$ induces a triangle-equivalence $D^b(B\mbox{-}{\rm
mod})/{\rm Ker} \pi_*\cong D^b(A\mbox{-}{\rm mod}).$

\vskip5pt

Notice that $\pi^*$ is left adjoint to $\pi_*$ on $K^b(\mathcal
P(A\mbox{-}{\rm mod}))$. In fact, for $P\in K^b(\mathcal
P(A\mbox{-}{\rm mod}))$ and $X\in D^{b}(B\mbox{-}{\rm mod})$ we have
\begin{align*}\Hom_{D^b(B\mbox{-}{\rm mod})}(\pi^*P, X)&\cong
\Hom_{K^{-, {\rm add M}\ b}({\rm add M})}(\sigma P, FX)\\ & \cong
\Hom_{K^{-, {\rm add M}\ b}({\rm add M})}(P, FX);\end{align*} and
\begin{align*}\Hom_{D^b(A\mbox{-}{\rm mod})}(P, \pi_*X) & \cong
\Hom_{K^{-, {\rm add M}\ b}({\rm add M})/K^{b}_{ac}({\rm add
M})}(GP, VFX)\\ & \cong \Hom_{K^{-, {\rm add M}\ b}({\rm add
M})/K^{b}_{ac}({\rm add M})}(P, FX)\end{align*} (note that $GP = P$
and $VF X = FX$). So, it suffices to prove that there is a
functorial isomorphism
$$\zeta_{P, FX}: \Hom_{K^{-, {\rm add M}\ b}({\rm add M})}(P, FX)\cong \Hom_{K^{-, {\rm add M}\ b}({\rm add M})/K^{b}_{ac}({\rm add M})}(P,
FX).$$ This follows from Proposition \ref{ff} by taking $\mathcal C
= {\rm add M}$.

\vskip5pt

Finally,  saying that the unit ${\rm Id}_{\mathcal D^{\rm perf}}
\longrightarrow \pi_*\pi^* = G^{-1} V \sigma$ is a natural
isomorphism of functors amounts to saying that
$$\zeta_P = \zeta_{P, P}({\rm Id}_P) = {\rm Id}_P/{\rm Id}_P: P\longrightarrow P$$ is an
isomorphism in $K^{-, {\rm add M}\ b}({\rm add M})/K^{b}_{ac}({\rm
add M})$ for each $P\in K^b(\mathcal P(A\mbox{-}{\rm mod}))$. This
is trivially true.

\vskip5pt

All together triple $(D^b(B\mbox{-}{\rm mod}), \pi_*, \pi^*)$ is a
categorical resolution of $D^b(A\mbox{-}{\rm mod}).$  $\s$

\vskip10pt

\subsection{} Theorem \ref{mainresult} is stated for the finitely
generated module category $A\mbox{-}{\rm mod}$. By the similar
argument with a minor change we can prove its version for
$A\mbox{-}{\rm Mod}$. For the contravariantly finiteness of ${\rm
Add} M$ in $A$-{\rm Mod}, we need the assumption that ``$M$ is
finitely generated" (cf. Example \ref{exm}$(i)$). For $T\in
A\mbox{-}{\rm mod}$, let $^{\perp_{\rm big}} ({\rm Add}T)$ denote
the full subcategory of $A$-Mod given by $^{\perp_{\rm big}} ({\rm
Add}T) = \{X\in A\mbox{-}{\rm Mod} \ | \ \Ext^i_A(X, T') = 0, \
\forall \ i\ge 1, \ \forall \ T'\in {\rm Add}T \}$ (here ``big"
refers to work in $A$-Mod). Note that there is no a module $T'\in
A$-Mod such that $^{\perp_{\rm big}} ({\rm Add}T) = \ ^{\perp_{\rm
big}} T'.$

\vskip5pt

\begin{thm}\label{finrem} Let $A$ be an Artin algebra with ${\rm gl.dim}A = \infty$. Assume
that there are modules $T$ and $M$ in $A\mbox{-}{\rm mod}$ with
${\rm inj.dim} T< \infty$, such that $^{\perp_{\rm big}} ({\rm
Add}T) = {\rm Add} M$.
 Then $D^b(A\mbox{-}{\rm Mod})$ admits a categorical
resolution $D^b(B\mbox{-}{\rm Mod})$ with $B = (\End_AM)^{op}$.
\end{thm}
\noindent{\bf Proof.} First, the condition $^{\perp_{\rm big}} ({\rm
Add}T) = {\rm Add} M$ implies $^\perp T = {\rm add}M$. The argument
is as follows:
$$^\perp T = \ ^{\perp_{\rm big}} ({\rm Add}T) \cap A\mbox{-}{\rm
mod} = {\rm Add} M\cap A\mbox{-}{\rm mod} = {\rm add}M.$$ By Theorem
\ref{global} ${\rm gl.dim} B < \infty,$ i.e., $D^b(B\mbox{-}{\rm
Mod})$ is smooth. Since $M$ is finitely generated, $\Hom_A(M, -):
{\rm Add M} \longrightarrow \mathcal P(B\mbox{-}{\rm Mod})$ is again
an equivalence of categories. Since
$$\Hom_A(X, M') \cong \Hom_B(\Hom_A(M, X), \Hom_A(M, M')), \ \forall \ X\in {\rm Add}M, \ \forall \  M' \in {\rm Add}M,$$
it follows that this equivalence induces pointwisely a
triangle-equivalence $K^{-, {\rm Add M}\ b}({\rm Add M})\cong K^{-,
b}(\mathcal P(B\mbox{-}{\rm Mod}))$, and hence we get a
triangle-equivalence
$$F: D^b(B\mbox{-}{\rm Mod}) \cong K^{-, {\rm Add M}\ b}({\rm Add M}).$$
Since ${\rm Add M} = \ ^{\perp_{\rm big}} ({\rm Add}T)$, it follows
that ${\rm Add M}$ is a resolving subcategory of $A$-Mod. Also,
${\rm Add M}$ is contravariantly finite in $A$-Mod by Example
\ref{exm}$(i)$.

\vskip5pt

The rest of the proof is similar with the one for Theorem
\ref{mainresult}, just replacing ${\rm add M}$ by ${\rm Add M}$,
$A\mbox{-}{\rm mod}$ by $A\mbox{-}{\rm Mod}$, and $B\mbox{-}{\rm
mod}$ by $B\mbox{-}{\rm Mod}$. We omit the details. $\s$

\vskip10pt

\subsection{} Let us see some
special cases of Theorems \ref{mainresult} and \ref{finrem}. We have
a reformulation of the Auslander algebra:

\vskip5pt

\begin{cor} \label{mainresult1} \  Let $A$ be a representation-finite Artin algebra with ${\rm gl.dim}A = \infty$, and
$B$ its Auslander algebra. Then

\vskip5pt

$(i)$ \ $D^b(B\mbox{-}{\rm mod})$ is a categorical resolution
$D^b(A\mbox{-}{\rm mod})$.

\vskip5pt

$(ii)$ \ $D^b(B\mbox{-}{\rm Mod})$ is a categorical resolution
$D^b(A\mbox{-}{\rm Mod})$.
\end{cor}
\noindent{\bf Proof.} Put $T$ to be an injective module in $A$-mod,
and $M$ to be the direct sum of all the pairwise non-isomorphic
finitely generated indecomposable modules.

By Theorem \ref{mainresult} we get $(i)$.

Since $A$ is representation-finite, any $A$-module is a direct sum
of finitely generated indecomposable modules (see [A2]). It follows
that $^{\perp_{\rm big}} T = A\mbox{-}{\rm Mod} = {\rm Add} M$. By
Theorem \ref{finrem} we get $(ii)$. $\s$

\subsection{} A module $T\in A$-mod is {\it a cotilting module} ([AR]), if
\vskip5pt

$(i)$ \ ${\rm inj.dim} T\le 1$;

\vskip5pt

$(ii)$ \ ${\rm Ext}^1_A(T, T) = 0$; and

\vskip5pt

$(iii)$ \ There is an exact sequence $0\rightarrow T_0 \rightarrow
T_1\rightarrow D(A_A) \rightarrow 0$ with $T_i\in {\rm add T}$, $i
=0, 1$.

\vskip5pt

An module $X\in A$-mod {\it is cogenerated by} $T$, if $X$ can be
embedded as an $A$-module into a finite direct sum of copies of $T$.
Then $X$ is cogenerated by a cotilting module $T$ if and only if
$X\in \ ^\perp T$ ([HR]). By Theorem \ref{mainresult} we have

\vskip5pt

\begin{cor} \label{mainresult2} \  Let $A$ be an Artin algebra with ${\rm gl.dim}A = \infty$. Assume that $A$ has a cotilting module $T$ such that
there are only finitely many pairwise non-isomorphic indecomposable
$A$-modules which are cogenerated by $T$. Then $D^b(A\mbox{-}{\rm
mod})$ admits a categorical resolution.
\end{cor}

\subsection{} Finally, we consider {\rm CM}-finite Gorenstein algebras.

\vskip5pt

\begin{thm} \label{mainresult3} \  Let $A$ be a {\rm CM}-finite Gorenstein algebra with ${\rm gl.dim}A =
\infty$, and $B$ its relative Auslander algebra. Then

\vskip5pt

$(i)$ \ $D^b(B\mbox{-}{\rm mod})$ is a weakly crepant categorical
resolution $D^b(A\mbox{-}{\rm mod})$.

\vskip5pt

$(ii)$ \ $D^b(B\mbox{-}{\rm Mod})$ is a weakly crepant categorical
resolution $D^b(A\mbox{-}{\rm Mod})$.
\end{thm}
\noindent {\bf Proof.}   Take $T = \ _AA$, and $M$ to be the direct
sum of all the pairwise non-isomorphic finitely generated
indecomposable Gorenstein-projective modules, in Theorem
\ref{mainresult} and \ref{finrem}.

\vskip5pt

$(i)$ \ Since $A$ is CM-finite, we have $M\in A\mbox{-}{\rm mod}$
and $\mathcal {GP}(A\mbox{-}{\rm mod}) = {\rm add} M$. Since $A$ is
Gorenstein, it follows from [EJ, Corollary 11.5.3] that $\mathcal
{GP}(A\mbox{-}{\rm mod}) = \ ^\perp (_AA)$. Thus $^\perp T = {\rm
add M}$. Then $D^b(A\mbox{-}{\rm mod})$ has a categorical resolution
$(D^b(B\mbox{-}{\rm mod}), \pi_*, \pi^*)$ by Theorem
\ref{mainresult}. It remains to see that $\pi^*$ is right adjoint to
$\pi_*$ on $K^b(\mathcal P(A\mbox{-}{\rm mod}))$. As in the proof of
Theorem \ref{mainresult} it suffices to prove that there is a
functorial isomorphism
$$\Hom_{K^{-, {\rm add M}\ b}({\rm add M})}(FX, P)\cong \Hom_{K^{-, {\rm add M}\ b}({\rm add M})/K^{b}_{ac}({\rm add M})}(FX, P).$$
This follows from Proposition \ref{inj3} by taking $\mathcal C =
{\rm add M} = \mathcal {GP}(A\mbox{-}{\rm mod})$, since projective
modules are injective objects of  $\mathcal {GP}(A\mbox{-}{\rm
mod})$.

\vskip5pt

$(ii)$ \ Since $A$ is a {\rm CM}-finite Gorenstein algebra, any
Gorenstein-projective $A$-module is a direct sum of finitely
generated indecomposable Gorenstein-projective modules (see [C1]).
It follows that $\mathcal {GP}(A\mbox{-}{\rm Mod}) = {\rm Add} M$.
Since $A$ is Gorenstein, it follows from  [EJ, Corollary 11.5.3] (or
[Bel1, Proposition 3.10]) that  $\mathcal {GP}(A\mbox{-}{\rm Mod}) =
\ ^{\perp_{\rm big}} ({\rm Add} \ _AA)$. Thus $^{\perp_{\rm big}}
({\rm Add}T) = {\rm Add} M.$  Then $D^b(B\mbox{-}{\rm Mod})$ is a
categorical resolution of $D^b(A\mbox{-}{\rm Mod})$ by Theorem
\ref{finrem}. By the similar argument as in $(i)$ we know that it is
weakly crepant. $\s$

\end{document}